\documentclass[12pt]{article}
\usepackage{amsmath, amsfonts, amssymb, amsthm} \parskip1pt
\newtheorem{theorem}{Theorem}[section]
\newtheorem{proposition}[theorem]{Proposition}
\newtheorem{definition}[theorem]{Definition}
\newtheorem{remark}[theorem]{Remark}

\newtheorem{corollary}[theorem]{Corollary}
\newtheorem{lemma}[theorem]{Lemma}

\def\R{\mathbb R} \def\Z{\mathbb Z}

\def\N{\mathbb N}
\def\C{{\mathbb C}} 
\def\Q{\mathbb Q}

\def\supp{{{\rm supp}\,}}
\def\wt{\widetilde}
\def\Ad{\hbox{\rm Ad}}

\def\wt{\widetilde}
\def\wh{\widehat}
\def\Hom{\hbox{\rm Hom }}

\begin{document}

\title{On the embeddability of certain infinitely divisible probability
  measures on Lie groups}

\author{S.G. Dani, Yves Guivarc'h and Riddhi Shah}
%{}\\}

\maketitle

\begin{abstract}
We describe certain sufficient conditions for an infinitely
divisible probability
measure on a class of  connected Lie groups to be embeddable in
continuous one-parameter semigroups of probability measures.
(Theorem~\ref{thm:main}). This
enables us in particular to conclude the embeddability of all
infinitely divisible
probability measures on certain Lie groups, including the so called
Walnut group (Corollary~\ref{walnut}).  The embeddability is
concluded also under certain other conditions (Corollary~\ref{G=NK}
and Theorem~\ref{noncomm}).

\end{abstract}

\section{Introduction}

Let $G$ be a Lie group. We denote by $P(G)$ the space of probability
measures on $G$, equipped with the convolution product and the usual
weak* topology. A $\mu\in
P(G)$ is said to be {\it infinitely divisible} if it admits
convolution roots of all orders and it is said to be {\it embeddable} if
there exists a continuous one-parameter convolution semigroup
$\{\mu_t\}_{t\geq 0}$ such that $\mu_1=\mu$. It has been conjectured
that when $G$ is
a connected Lie group every infinitely divisible probability measure
on $G$ is embeddable; the latter is known as the embedding property.
The results in \cite{DM-inv} and \cite{DM-adv} show that this holds
under the additional
condition that the nilradical of $G$ contains no torus of
positive dimension. The embedding property is also known to hold
for all nilpotent Lie groups. An example of a nonnilpotent Lie
group for which the above condition fails occurs in dimension~4.
Partial results towards the embedding property were proved for this
group,  termed as the Walnut group, in
\cite{DMW}, taking off from the Ph.D. thesis of S.~Walker.

In this paper we prove some new results on the embeddability of
infinitely divisible
measures on a class of  connected Lie groups, under an assumption
on the measure (see Theorem~\ref{thm:main}). This in particular
enables us to settle the embedding problem for the Walnut group, and
also a larger group containing it
(see Corollary~\ref{walnut}).

We now begin with some more notation and a formal statement of the
main results. Let $G$ be a Lie group.
For any closed subset $S$ of $G$ we denote by
$P(S)$ the subspace of $P(G)$ consisting of the probability measures
whose support is contained in $S$.
For a subset $\Lambda$ of probability measures on $G$  we denote by
$G(\Lambda)$ the smallest
closed subgroup of $G$ containing the supports of all $\lambda \in
\Lambda$, and by $Z(\Lambda)$  the centraliser of $G(\Lambda)$ in $G$.
For $\lambda \in P(G)$, $G(\{\lambda\})$ and $Z(\{\lambda\})$ will
be written as $G(\lambda)$ and $Z(\lambda)$ respectively.
For $\lambda \in P(G)$ we denote by
$N(\lambda)$ the normaliser of $G(\lambda)$ in $G$.
 We recall that given  a (convolution) root $\rho $ of $\lambda$,
 there exists $x\in N(\lambda)$ such that $\rho
\in P(xG(\lambda))$; furthermore, if $\rho^n=\lambda$ then for every such
$x$ we have $x^n\in G(\lambda)$ (see, for instance, \cite{DS},
Lemma~2.2).

We shall mainly be concerned  with Lie groups which admit
surjective homomorphisms onto almost algebraic groups, with kernels
contained in the center; see Theorem~\ref{thm:main} below.
By an almost algebraic group
we mean a Lie group which is (Lie isomorphic to) a subgroup of finite
index in  a real algebraic group (the group of real points of an
algebraic group defined over $\R$); see
\cite{DM-inv}, \S~2, for some generalities in this respect.
Let $\wt G$ be an almost algebraic group.
 For a subset $\Lambda$ of $P(\wt G)$ we shall denote by $\wt
G(\Lambda)$  the smallest almost algebraic subgroup
of $\wt G$ containing $G(\Lambda)$ (or, equivalently, $\supp
\lambda$ for all $\lambda \in \Lambda$). For  $\nu\in P(\wt G)$,
$\wt G(\nu)$ stands for  $\wt G(\{\nu\})$ and
$\wt N(\nu)$  denotes the normaliser of $\wt G(\nu)$ in $\wt G$.
We note that for any $\nu\in P(\wt G)$, $Z(\nu)$ is an algebraic
subgroup of $\wt G$ and it is also the centraliser of $\wt
G(\nu)$. Also, for any $\nu\in P(\wt G)$, $N(\nu)$ is contained in
$ \wt N(\nu)$.

For a Lie subgroup $H$ we denote by $H^0$ the connected component of
the identity in $H$; (in the case of $G(\lambda)$, $Z(\lambda)$ etc. the
identity component will be written as $G^0(\lambda)$, $Z^0(\lambda)$ etc.).

\begin{definition}
{\rm Let $\mu \in P(G)$. A subset $\Psi$ consisting of roots of $\mu$
  in $P(G)$ is called a {\it root cluster} of $\mu$ if the following
  conditions are satisfied:

i) if $n\in \N$ and $\rho \in \Psi$ is a $n$ th root of $\mu$ then
$\rho^k\in \Psi$ for all $k$ dividing $n$;

ii) for all $\rho \in \Psi$ and $z\in Z(\mu)$, $z\rho z^{-1} \in
\Psi$;

iii) for all $n\in \N$, $\{\rho \in \Psi \mid \rho^n=\mu\}$ is a
nonempty closed subset of $P(G)$.
}
\end{definition}

We note that if $\mu $ is infinitely divisible then the set of all
roots of $\mu$ in $P(G)$ is a root cluster for $P(G)$.

\begin{definition}
{\rm Let $\mu \in P(G)$. A family $\{\mu_r\}_{r\in \Q^+}$, where
$\Q^+$ is the additive semigroup consisting of all positive rationals,
is called a
  rational embedding of $\mu$ if $\mu_1=\mu$ and $r\mapsto \mu_r$ is a
  homomorphism of $\Q^+$ into $P(G)$.
}
\end{definition}

It is known that every rational embedding  $\{\mu_r\}_{r\in \Q^+}$
extends uniquely to an embedding  $\{\mu_t\}_{t\geq 0}$ (see
\cite{DM-jtp}), so if $\mu$ admits a rational embedding then it is
continuously embeddable.

We prove the following.

\begin{theorem}\label{thm:main}
Let $G$ be a connected Lie group admitting a surjective
continuous homomorphism $p:G\to \wt G$ onto an almost algebraic group
$\wt G$, such that $\ker p$ is contained in the center of $G$ and
$(\ker p)^0$ is compact. Let $T=(\ker p)^0 $  and $q:G\to
G/T$ be the quotient homomorphism. Let $\mu\in P(G)$ be such that
$q(\mu)$ has no nontrivial idempotent factor in $P(G/T)$. Let $\Psi$
be a root cluster of $\mu$. Then there exists
a  rational embedding $\{\mu_r\}_{r\in \Q^+}$ of
 $\mu$ such that $\mu_r\in \Psi$ for all $r\in \Q^+$; in particular, if
 $\mu$ is infinitely divisible then it
 is embeddable.
\end{theorem}

From the theorem we deduce also the following corollary,
for a general connected Lie group.

%\medskip
\begin{corollary}\label{G=NK}
Let $G$ be a connected Lie group and $N$ be the nilradical of
$G$.  If $\mu \in P(N)$ is infinitely divisible in $P(G)$ then $\mu$ is
embeddable in $P(G)$.
\end{corollary}

We next describe certain special groups for which the embedding problem
can now be settled.   Let $H$ be the
 3-dimensional Heisenberg group, $Z$ be the (one-dimensional) center
 of $H$  and $D$ an infinite cyclic subgroup of $Z$. Let $N=H/D$ and
$T=Z/D$. Then $N$ is a nilpotent Lie group and $T$ is the maximal
torus in $N$.
The group $SL(2,\R)$ has a canonical action on $H$ as a group of
automorphisms, such that $Z$ is point-wise fixed and the quotient
action on $N/Z$, which we view as $\R^2$, is the usual linear
action of $SL(2,\R)$ on $\R^2$. As $D$ is pointwise fixed under the
action we get a quotient action of $SL(2,\R)$ on $N$, as a group
of automorphisms. Let $G$ be the semidirect product of $SL(2,\R)$ and
$N$ with respect to the action. It may be recalled here that this
group plays an
important role in the study of Jacobi forms; see \cite{EZ}.
We view $SL(2,\R)$ and $N$ as subgroups of $G$,
canonically. Let $K$ be the subgroup of $SL(2,\R)$ consisting of
rotations of $\R^2$ around $0$ (with respect to a fixed Hilbert norm
on $\R^2$).
The subgroup of $G$ containing $K$ and $N$, which is the semidirect
product of the two subgroups with respect to the $K$ action on $N$,
is the {\it Walnut} group; see \cite{DMW} and \cite{M-cimpa} for more
details. We shall denote the subgroup  by $W$.

Applying Theorem~\ref{thm:main} together with certain known results we
deduce the following.

\begin{corollary}\label{walnut}
Let $G$ be the semidirect product of $SL(2,\R)$ and $N$ as above.
Let $H$ be a closed connected subgroup of $G$. Then
every infinitely divisible probability measure in  $P(H)$ is
embeddable in  $P(H)$. In particular, every probability measure which
is infinitely divisible in  $P(W)$ is embeddable in $P(W)$.
\end{corollary}

In the proof of Theorem~\ref{thm:main} the condition that $q(\mu)$ has no
idempotent factor is needed only in the case
$G(\mu)/\overline{[G(\mu),G(\mu)]}$ is not compact. Using this
observation together with some facts about Lie groups we conclude also
the following result, for all connected Lie groups (not involving the
condition as in Theorem~\ref{thm:main}).

\begin{theorem}\label{noncomm}
Let $G$ be a connected Lie group and $\mu \in P(G)$ be infinitely
divisible. Suppose that $G(\mu)/\overline{[G(\mu),G(\mu)]}$ is
compact. Then $\mu$ is embeddable.
\end{theorem}

\section{Some generalities on central  extensions}

Let $G$ be a Lie group, not necessarily connected, and suppose that
there exist an
almost algebraic group $\wt G$ and a surjective continuous
homomorphism $p:G\to \wt G$ such that $\ker p$ is a
compactly generated subgroup contained in the
center of $G$ and $(\ker p)^0$ is compact. Then $\ker p$ has a
unique maximal compact subgroup which we shall denote by $A$, and
$\ker p$ is isomorphic to $A \times \Z^r$ for some $r\geq 0$. The
connected component of the identity in $A$, which is a torus, will be
denoted by $T$.

 Let $\mu \in P(G)$ and $\nu=p(\mu)$.
For each $z\in Z(\nu)$ we define a map  $\psi_z:\wt
G(\nu) \to \ker p$ as follows: let  $x\in \wt G(\nu)$, and $\xi, \zeta \in G$
be such that $p(\xi)=x$ and $p(\zeta)=z$ (we recall that $p$ is surjective);
since $x$ and $z$ commute with each other, the element
$\zeta \xi \zeta^{-1}\xi^{-1}$ is contained in $\ker p$, and since $\ker p$
is central in $G$ it follows that this element is independent
of the choices of $\xi$ and  $\zeta$  (depends only on $x$ and $z$);
we define  $\psi_z(x)$ to be  $\zeta \xi \zeta^{-1}\xi^{-1}$.  Clearly
$\psi_z$ is a continuous map for all $z$. Also, for $z=p(\zeta)\in
Z(\nu)$, given
$x=p(\xi), x'=p(\xi')\in \wt G(\nu)$, we have $\psi_z(xx')=\zeta \xi \xi'
\zeta^{-1}
(\xi\xi')^{-1}=\zeta \xi \zeta^{-1}\psi_z(x')\xi^{-1}=\psi_z(x)\psi_z(x'),$
as $\psi_z(x')$ is central. This shows that each $\psi_z$  is a
homomorphism. Being an almost algebraic group  $\wt G (\nu)$ has
only finitely many connected components, and hence it follows also
that for all $z \in
Z(\nu)$ the image of $\psi_z$ is contained in $A$. Furthermore, as the
image is abelian, $\psi_z$ factors canonically to a  homomorphism of $\wt
G(\nu)/[\wt G(\nu), \wt G(\nu)]$ into  $A$ (we note that as
$\wt G(\nu)$ is almost algebraic, $[\wt
G(\nu), \wt G(\nu)]$ is also almost algebraic and hence in particular
a closed subgroup);  we shall denote this homomorphism also
by $\psi_z$.

\begin{lemma}\label{psi}
Let $Q$ be the smallest closed normal subgroup of $\wt G(\nu)$ such
that $\wt G(\nu)/Q$ is a vector group, and let $\wt V=\wt
G(\nu)/Q$. Let $K$ be the
smallest closed normal subgroup of $Z^0(\nu)$ such that
$Z^0(\nu)/K$ is a vector group. Then we have the following.

i) for all $z\in Z^0(\nu)$, $ \psi_z$ factors to a continuous homomorphism
of $\wt V$ into $T$;

ii) for all $z\in K$, $\psi_z$ is trivial, that is $\psi_z(x)=e$ for
all $x\in \wt G(\nu)$.

\end{lemma}

\proof Let $\Hom
(\wt G(\nu)/[\wt G(\nu), \wt G(\nu)], A)$ be the group of
continuous homomorphisms of $\wt G(\nu)/[\wt G(\nu), \wt G(\nu)]$ into
$A$, under pointwise multiplication, equipped
with the compact open topology. Let  $\psi : Z(\nu) \to \Hom (\wt G(\nu)/[\wt
G(\nu), \wt G(\nu)], A)$ be the map $z\mapsto \psi_z$ for all $z\in \wt
Z(\nu)$ defined as above. It can be seen that $\psi$ is a continuous
homomorphism.

As $\wt G(\nu)$ is an almost algebraic group,  $\wt G(\nu)/[\wt
G(\nu), \wt G(\nu)] $ is an abelian Lie group with finitely many
connected components. By Pontryagin duality theory this implies that
$\Hom (\wt G(\nu)/[\wt G(\nu), \wt G(\nu)], A)$
has no compact subgroup of positive dimension. In particular its
connected component of the identity is a vector group, (namely a
topological group isomorphic to $\R^n$ for some $n\geq 0$). The identity
component therefore consists of homomorphisms that factor to $\wt V$
and have their image
contained in $T$; we shall denote the latter subgroup by $\Hom (\wt V, T)$.
Then by continuity, $\psi (Z^0(\nu))$ is contained in $\Hom (\wt V, T)$;
that is, for all $z\in Z^0(\nu)$, $ \psi_z$ factors to a continuous
homomorphism of $\wt V$  into $T$. This proves (i).
Assertion~(ii) is immediate from the definition of $K$ and the fact
that $\psi (Z^0(\nu))$ is a vector group. This proves the Lemma. \hfill $\Box$

\begin{remark}\label{invrem}
{\rm If $\zeta \in p^{-1}(Z(\nu))$ is such that
$\psi_{p(\zeta)}$ is trivial then $\zeta \in Z(\mu)$. Indeed, for any $\xi \in
G(\mu)$  we have $\zeta \xi \zeta^{-1}\xi^{-1}= \psi_{p(\zeta)} (p(\xi))=e$.
Lemma~\ref{psi} implies that $(z,x)\mapsto \psi (z,x)$ is a continuous
bi-homomorphism of $\wt V \times
Z^0(\nu)/K$ into $T$. Using Pontryagin duality we see that this map
defines a bilinear pairing into the universal cover of $T$.
}
\end{remark}

Let $K$ be the subgroup as in Lemma~\ref{psi}. We note that $K$ is
the subgroup containing $[Z^0(\nu),Z^0(\nu)]$ (which is an almost
algebraic normal subgroup of $Z^0(\nu)$) and such that
$K/[Z^0(\nu),Z^0(\nu)]$ is the maximal torus in the abelian Lie
group $Z^0(\nu)/[Z^0(\nu),Z^0(\nu)]$. Let $W$ denote the vector space
$Z^0(\nu)/K$.
For $x\in \wt N(\nu)$ let $\tau_x:W\to W$ be the linear transformation
of $W$ induced by the
conjugation action of $x$ on $Z^0(\nu)$, and let $W_x$ be the range of
$\tau_x-I$, where $I$ denotes the identity transformation of $W$. We
note that
$W_x$ is the smallest $\tau_x$-invariant subspace such that the factor
of $\tau_x$ on $W/W_x$ is the identity automorphism.
We shall denote by $K(\nu, x)$ the  subgroup of $Z^0(\nu)$
containing $K$ and such that $W_x=K(\nu,x)/K$.

\begin{lemma}\label{lem:comm}
Let $x\in \wt N(\nu)$ and $\{g_t\}$ be a one-parameter subgroup
in $\wt G(\nu)$ such that $xg_tx^{-1}g_t^{-1}$ is contained in
$Q$, namely $\{g_tQ\}$ which is a one-dimensional subspace of $\wt
V=\wt G(\nu)/Q$ is invariant under the action of $x$ on $\wt V$
induced by the conjugation action. Then $\psi_z(g_t)=e$ for all
$t\in \R$ and $z\in K(\nu, x)$.
\end{lemma}

\proof Let $\mathfrak G$ denote the Lie algebra of $G$. Let $\xi \in
\mathfrak G$ be  such that $p(\exp t\xi)=g_t$ for all $t\in \R$ and let
$\zeta \in \mathfrak G$ be any element such that $\{p(\exp t\zeta)\}$ is
contained in $Z^0(\nu)$. Then we have $\Ad \, x
([\xi,\zeta])=[\Ad \,x
(\xi),\Ad x\, (\zeta)]= [\xi+\theta,\Ad\, x (\zeta)]$, where $\theta $
is an element of the Lie subalgebra of $p^{-1}(Q)$. Since the
restriction of $\psi_z$ to $Q$ is trivial for all $z\in Z^0(\nu)$ (by
Lemma~\ref{psi}) it
follows that $[\theta,\Ad\, x (\zeta)]=0$, and hence we have  $\Ad \, x
([\xi,\zeta])= [\xi,\Ad\, x (\zeta)]$. On the other
hand since $[\xi,\zeta]$ belongs to the Lie subalgebra corresponding
to $T$, which is contained in the center of $G$, we have $\Ad \, x
([\xi,\zeta])= [\xi,\zeta]$. Therefore $[\xi,\Ad \, x (\zeta)]=
[\xi,\zeta] $, or equivalently $[\xi ,\Ad \, x (\zeta)-\zeta]=0$ for
all $\zeta$ as above.
This shows that for all $t\in \R$, $\exp t\xi$ centralises
$K(\nu, x)$. Equivalently $\psi_z(g_t)=e$ for all
$t\in \R$ and $z\in K(\nu, x)$. \hfill $\Box$

\medskip
We need the following consequence of the lemma in the proof of
Theorem~\ref{thm:main}.

\begin{corollary}\label{cor:nofp}
Let $\Phi$ be a set of roots of $\nu$. Let $M$ be the smallest
closed normal $\wt G (\Phi)$-invariant subgroup of $Z^0(\nu)$
 such that $Z^0(\nu)/M$ is a vector group and
the induced $\wt G(\Phi)$-action on $ Z^0(\nu)/M$ is trivial. Let $g\in
\wt G(\nu)$ be such that $gQ$ is a fixed point of the
$\wt G (\Phi)$-action. Then $\psi_z(g)=e$ for all $z\in M$.
\end{corollary}

\proof Since $\wt G(\nu)/Q$ is a vector group, the subspace spanned by $gQ$
is pointwise fixed under the $\wt G(\Phi)$-action and hence there
exists a one-parameter subgroup $\{g_t\}$ in $\wt G(\nu)$ such that
$g_1=g$ and $g_tQ$ is pointwise fixed by $\wt G(\Phi)$ for all $t\in
\R$. Then for any $x\in \wt G(\Phi)$,  by Lemma~\ref{lem:comm}
we get that $\psi_z(g_t)=e$ for all $t\in \R$ and $z\in
K(\nu, x)$, in the notation as above. Let $W'$ be the subspace spanned
by all $W_x$, $x\in \wt
G(\Phi)$. Then the preceding conclusion (together with the definition
of $W_x$) implies that $\psi_z(g)=e$ for all $z\in Z^0(\nu)$ such that
$zK\in W'$. Now, $W'$ is invariant under the $\wt G(\Phi)$ action, and
since $W_x$ is contained in $W'$ for all $x\in \wt G(\Phi)$ it follows
that the $\wt G(\Phi)$-action on the quotient $W/W'$ is trivial. Hence
$M$ as in the hypothesis is contained in $W'$. This proves the
corollary. \hfill $\Box$

\section{Measures on vector spaces}

In this section we note certain properties of measures on vector
spaces. For an abelian Lie group $H$ we denote by $\wh H$ its dual group and
for a finite measure $\lambda$ on $H$ we
denote by $\widehat
\lambda$ the Fourier transform of $\lambda$ defined on $\wh H$. A
function $f$ on $\widehat
H$ is said to {\it vanish at infinity} if for all $\epsilon >0$ there exists
a compact subset $C$ of $\widehat H$ such that $|f(\chi)|<\epsilon $ for
all $\chi \notin C$.

\begin{proposition}\label{lem:FT}
Let $V =\R^d$, $d\geq 2$, and $C$ be a compact connected Lie
subgroup of $GL(V)$ whose action has no nonzero fixed point in $V$.
Let $\lambda \in P(V)$ be a measure invariant under the
$C$-action.  Let $A$ be an affine subspace of $V$ such that $\lambda
(A)>0$.  Then $A$ is a vector subspace of $V$.

\end{proposition}

\proof  Let $U$
be a minimal affine subspace contained in $A$ with $ \lambda (U)>0$.
By the minimality condition for any $g\in C$
 either $gU=U$ or $\lambda (U\cap gU)=0$. Let $C'=\{g\in C\mid
 gU=U\}$; then $C'$ is a closed subgroup of $C$. Suppose that $C'$ is
 a proper closed subgroup. As $C$ is a connected Lie group it follows
that $C/C'$ is uncountable. We can therefore find an uncountable
subset $E$ of $C$ consisting of elements belonging to distinct cosets
of $C'$, so for any distinct $g,g'$ in $E$ we have $g^{-1}g'\notin
C'$. Then  $\lambda (gU\cap g'U)=\lambda (U\cap g^{-1}g'U) =0$ for
any distinct $g,g'$ in
 $E$, which means that $\{gU\}_{g\in E}$ are pairwise essentially
disjoint.  But this is not possible since $\lambda$ is $C$-invariant,
 $\lambda (U)>0$ and $E$ is uncountable. Therefore $C'=C$, which means
 that $U$ is $C$-invariant.
Now we consider $V$ equipped with a $C$-invariant Hilbert norm. Then
the point of $U$ closest to the origin would be a fixed point for the
$C$-action. Since by hypothesis there is no nonzero fixed point, the point
must be the origin.  Therefore $U$, and in turn $A$, contains $0$.
This shows that $A$ is a vector subspace.

\begin{theorem}\label{yves}
Let $V =\R^d$, $d\geq 2$, and $C$ be a compact connected Lie
subgroup of $GL(V)$ whose action has no nonzero fixed point in $V$.
Let $\lambda \in P(V)$ be a $C$-invariant measure on $V$, such that
there is no proper vector subspace of positive $\lambda$-measure. Then
$\lambda^d $ has a density on $V$ with respect to the Lebesgue measure;
in particular, $\widehat \lambda$ vanishes at infinity.

\end{theorem}

\proof For all $x\in V\setminus \{0\}$, $C x$ is an analytic
submanifold of $V$ of positive dimension. For any $x\in V$ we denote by
$T_{x}$ the tangent space
 to $C x$ at $x$; then  dim $T_{x}>0$ for all $x\neq 0$. We shall
first show that for almost all $(x_1,\dots , x_d)$ in $V^d$, with
respect to the measure $\lambda^{\otimes d}$ (the $d$th Cartesian
power of $\lambda$),  $T_{x_{1}}+\cdots + T_{x_{d}}=V$. We observe that, if
$x_{1}, x_{2},\cdots , x_{n}$ is a
finite sequence in $V$ such that $T_{x_{k+1}}$ is not contained in
$\sum_{i=1}^k\, T_{x_{i}}$ for $k= 1,\dots,  n-1$, then
$k\mapsto \dim (\sum _{i=1}^k\, T_{x_{i}})$ is strictly
increasing, hence $n\leq d$. For $\omega=(x_{1},\cdots , x_{d+1}) \in
V^{d+1}$ we denote by $n(\omega)$ the least $k\leq d$ such that
$T_{x_{k+1}}\subset \sum _{i=1}^k\, T_{x_{i}}$, and by  $p(\omega)$
the dimension of $\sum _{i=1}^{n(\omega)}\, T_{x_{i}}$; then
$n(\omega)\leq p(\omega)\leq d$. Let $p$ be the essential infimum of
the function $p(\omega)$ over $V^{d+1}$ with respect to the measure
$\lambda^{\otimes (d+1)}$. If we show that $p=d$ it would follow
that  $\sum _{i=1}^d\, T_{x_{i}}=V$, $\lambda^{\otimes d}$-a.e..

Suppose $p<d$. Then there exist $x_1,\dots , x_k$, for some $k\leq p$
such that $ \sum _{i=1}^k\, T_{x_{i}}$ is a $p$-dimensional subspace
and $T_x \subset  \sum _{i=1}^k\, T_{x_{i}}$ for all $x$ in a set of
positive $\lambda$ measure. Let $W= \sum _{i=1}^k\, T_{x_{i}}$ and
$S=\{x\in V\mid T_x\subset W\}$. Then we have $\lambda
(S)>0$. For $x\in V$ let $m_{\bar{x}}$ denote the normalised
$C$-invariant measure on $Cx$, let $\bar{V}$ be the space
of $C$-orbits and  $\bar{\lambda}$ the projection of $\lambda$
on $\bar{V}$. We decompose $\lambda$ into conditional measures
over the quotient space. Since $\lambda $ is $C$-invariant the conditional
measures are $ m_{\bar{x}}$ for $\bar \lambda$-almost all $\bar x\in
\bar V$. Thus
$\lambda=\int m_{\bar{x}}\,  d\bar{\lambda}(\bar{x})$.
Let $S'=\{x\in S \mid m_{\bar{x}} (S) >0\};$ then $\lambda
(S')=\lambda (S)>0$. Consider $x\in S'$. Since
the map $y\mapsto T_{y}$ is analytic over $y \in Cx$, and $T_y\subset W$ for
$y$ in a set of positive $ m_{\bar{x}}$-measure it follows that
$T_{y}\subset W$ for all $y\in Cx$. Now if
$\{\gamma (t)\}$ is a  differentiable
curve on $Cx$ and $l$ is a linear form on $V$ vanishing on
$W$, then the condition that $T_y\subset W$ for all $y\in Cx$ implies
that the
derivative of $l \circ \gamma$ is zero. This shows that $Cx$
is contained in $x+W$. Let $U$ be the subspace spanned by
$\{kx-k'x\mid k,k'\in C\}$. Then $U$ is a $C$-invariant subspace of
$W$ and also $x+U$ is $C$-invariant. We consider $V$ equipped with a
$C$-invariant Hilbert norm. Then $x+U$ has a unique point of minimum
norm, and hence the point must be fixed under the $C$-action. Since
by hypothesis there are no nonzero fixed points we get that $0\in
x+U$, which means that $x\in U\subset W$. Thus we have shown that
$S'$ is contained in $W$. Hence $\lambda (W)\geq \lambda
(S')>0$. Since $W$ is of dimension $p<d$ and by hypothesis there is no
proper subspace of positive
$\lambda$-measure this is a contradiction. This implies
that $p=d$, and hence
$  \sum _{i=1}^k\, T_{x_{i}}=V$, $\lambda^{\otimes d}$ - a.e..

We now consider the map $\varphi:V^d \to V$ given by
$\varphi(\omega)=x_{1}+\cdots+x_{d}$ for $\omega=(x_{1}, \cdots ,
x_{d}) \in V^d$, and the Borel set $E=\{\omega=(x_{1}, \cdots ,
x_{d}) \in V^d\mid  \sum_{i=1}^{d} T_{x_{i}}=V\}$ which satisfies
$\lambda^{\otimes d}(E)=1$. For $\omega \in E$, the differential  of
$\varphi$ at $\omega$ is surjective, and  hence for some neighbourhood
$\Delta$ of $\omega$ in $C{x_{1}}\times \cdots \times C{x_{d}}$,
$\varphi (\Delta)$ is an open subset of $V$; also the push forward by
$\varphi$ of the restriction to $\Delta$ of $m_{\bar{x}_{1}} \otimes
\cdots \otimes m_{\bar{x}_{d}}$ has a density with respect to Lebesgue
measure on $V$. Since this property is valid $\lambda^{\otimes d}$ -
a.e., and $\lambda^d=\int \varphi (m_{\bar{x}_{1}}\otimes \cdots
\otimes m_{\bar{x}_{d}}) d \lambda^{\otimes d} (\omega)$, we get that
$\lambda^d$ has also a density with respect to Lebesgue measure on
$V$. The last assertion in the statement of the theorem follows from
the Riemann-Lebesgue Lemma.  \hfill $\Box$

\section{Asymptotics of measures under shear transformations}

In this section we apply Proposition~\ref{lem:FT} to prove certain
general results about asymptotics of measures under the
action of sequences from a special class of automorphisms; see
\cite{D-cimpa} for other results on the theme.

Let $L$ be a locally compact
second countable group. Suppose that $L$ has a torus $T$ contained in
its center, and a closed normal subgroup $R$ such that $L/R$ is a
vector group. Let $W=L/R$ and $\theta :L\to W$ be the
quotient homomorphism. For any $\alpha \in \Hom (W,T)$ we get an
automorphism
$\tilde \alpha$ of $L$ defined by $\tilde \alpha (g)=g\alpha (gR)$ for
all $g\in G$; such an automorphism is  called a {\it shear automorphism}.
%Let $\Gamma$ be the group of automorphisms $\{\tilde \alpha \mid \alpha \in
%\Hom (W,T)\}$ and

We consider below $\Hom (W,T)$ as a vector space; it will be
convenient to use the additive notation, with the symbol ``$+$'' for the group
operation in $\Hom (W,T)$. For any $\mu\in P(L)$ let
$X(\mu)=\{ \alpha \in \Hom (W,T)\mid \tilde \alpha (\mu)=\mu\}$.

\begin{lemma}\label{sub}
Let the notation be as above. Let $\mu \in P(L)$. Suppose that every
affine subspace of
$W$ of positive $\theta (\mu)$-measure is a vector subspace. Then
$X (\mu)$ is a vector subspace of $\Hom (W,T)$.
\end{lemma}

\proof Let  $\alpha \in X (\mu)$ be a nontrivial element and
let $\{\alpha_t\}$ be the unique line in $\Hom (W,T)$ containing
$\alpha$. Let $\{\mu_w\}_{w\in W}$ be the conditional measures of
  $\mu$ over the fibration $\theta :L \to W$. We note that for
any $\alpha \in X(\mu)$ and $w\in W$, $\theta^{-1}(w)$ is
$\alpha$-invariant. It follows that for   $\theta (\mu)$-almost all
$w$, $\mu_w$
  is invariant under translation by $\alpha (w)$. The set of $w$ in
  $W$ for which the subgroup generated by $\alpha (w)$ is not dense in
  the closed subgroup generated by $\{\alpha_t(w)\mid t\in \R\}$ is
  a countable union of affine subspaces which are not vector
  subspaces; to see this note that for such a $w$ there exists a
  character $\chi$ on $T$ such that $\chi (\alpha (w))$ is a
  nontrivial element of finite order. Since   affine subspaces other
  than vector subspaces are
of $\theta (\mu)$-measure~$0$ this implies that for $\theta
(\mu)$-almost all $w$, $\mu_w$ is invariant under translation by
$\{\alpha_t(w)\}$. Hence $\mu$ is $\{\tilde
\alpha_t\}$-invariant, so $\alpha_t \in X (\mu)$ for all $t\in
\R$. The conclusion shows that $X (\mu)$ is connected, so being a
closed subgroup of $\Hom (W,T)$ it is a vector subspace.~\hfill $\Box$

\medskip
 A measure $\lambda$ on $V=\R^d$, $d\geq 2$ is said to
be {\it pure} if every vector subspace with positive $\lambda$-measure
contains the subspace spanned by $\supp \lambda$;  the
latter can be a proper subspace. It can be seen that
every $\lambda \in P(V)$ can be expressed uniquely as a countable sum
$\lambda =\sum_{i\in I} \lambda_i$ where $\{\lambda_i\}_{i\in I}$ is a
family, finite or countably infinite, of pure measures on $V$, and for
any distinct $i,j\in I$ the subspaces spanned by $\supp \lambda_i$
and $\supp \lambda_j$ are distinct; we call
 $\lambda_i$'s the {\it pure components} of~$\lambda$. For ready
 reference later we note the following.
%If $\lambda$ is invariant under the
%action of a connected Lie subgroup of  $GL(V)$ then each pure component
%of $\lambda$ is also invariant under the action of the subgroup.
\iffalse
\begin{remark}\label{inv-comp}
{\rm Let $\lambda \in P(\R^n)$, $n\geq 2$, and let $\lambda
  =\sum_{i\in I} \lambda_i$, where $I$ is an indexing set, be the
  decomposition of $\lambda$ into
  pure components. Let $C$ be a connected subgroup of $GL(n,\R)$, and
suppose that $\lambda$ is invariant under the action of $C$ on
$\R^n$. Then $\lambda_i$ is invariant under the $C$-action for all
$i\in I$. We note that by uniqueness of the decomposition the action
of any $g\in C$ must permute the components $\lambda_i$, $i\in
I$, and furthermore, the permutation must leave invariant, for every
$a>0$, the finite set of $\lambda_i$'s for which $\lambda_i(\R^n)=a$.
As the group is connected this implies that the
permutation must be trivial.
}
\end{remark}
\fi
\begin{lemma}\label{inv-comp}
Let $\lambda \in P(\R^n)$, $n\geq 2$, and let $\lambda
  =\sum_{i\in I} \lambda_i$ be the decomposition of $\lambda$ into
  pure components, where $I$ is an indexing set. For all $i$ let $W_i$
  be the subspace spanned by $\supp \lambda_i$. Then the following
  statements hold:

i) If $\lambda$ is invariant under the action of a connected subgroup
$C$ of $GL(n,\R)$ then each $\lambda_i$ is $C$-invariant.

ii) If for all $i$ there exist $d_i\in \N$ such that $\lambda_i^{d_i}$
has a density with respect to the Lebesgue measure on $W_i$, then
every affine subspace of $\R^n$ with positive $\lambda$-measure is a
vector subspace.
\end{lemma}

\proof i) Let $C$ be a connected subgroup of $GL(n,\R)$ such that
$\lambda$ is $C$-invariant. By the uniqueness of the collection of
pure components $\lambda_i$, $i\in I$, $C$ must permute the
components. For any $i\in I$ the subgroup $C_i$ of elements of
$C$ leaving $\lambda_i$ invariant is a closed subgroup, and as $I$ is
countable we get that $C/C_i$ is countable; since $C$ is connected
this implies that $C=C_i$, so $\lambda_i$ is $C$-invariant.

ii) Let $A$  be an affine subspace such that $\lambda
(A)>0$. Let  $i\in I$ be such that $\lambda_i (A)>0$. Let
$B=\{\sum_{j=1}^{d_i} v_j\mid v_j\in
A\}=\{d_iv\mid v\in A\}$. Then $\lambda_i^{d_i}(B)\geq
\lambda_i(A)^{d_i}>0$. Since $\lambda_i^{d_i}$ is supported on $W_i$
and has a density on $W_i$ this implies $B$ contains $W_i$. In
particular $0\in B$, and since $B=d_iA$ we get that $0\in
A$. Therefore $A$ is a vector subspace.

\begin{theorem}\label{thm:FT}
 Let $\mu$ be a finite measure on $L$.
%such that all affine subspaces of $V$ with positive $\theta
%(\mu)$-measure are  vector subspaces.
Let $\theta (\mu) =\sum_{i\in
  I} \lambda_i$ be the decomposition of $\theta (\mu)$ into pure components
and for
all $i$ let $W_i$ be the subspace spanned by $\supp \lambda_i$. Suppose that
for all $i$ there exists $d_i$ such that $\lambda_i^{d_i}$
has a density with respect to the Lebesgue measure on
$W_i$.  Let $\{\alpha_j\}$ be a
sequence in $\Hom (W,T)$ such that $\{\tilde \alpha_j(\mu)\}$ converges to
a measure of the form $\tau (\mu)$, for some continuous automorphism
$\tau$ of $L$ fixing  $T$ pointwise. Then in the quotient space $\Hom
(W,T)/X (\mu)$ the sequence  $\{\alpha_jX(\mu)\}$  is
relatively compact.
\end{theorem}

\proof The measure $\mu$ can be expressed uniquely as  $\mu=\sum_{i\in I}
\mu_i$, where  $\mu_i\in P(L)$, $i\in I$, are such that $\theta
(\mu_i)=\lambda_i$ for all $i$. Then $X (\mu)
=\cap X (\mu_i)$. By Lemmas~\ref{sub} and \ref{inv-comp} the condition
in the hypothesis implies that each $X (\mu_i)$ is a vector
subspace of $\Hom (W,T)$. Therefore to show that
$\{\alpha_j+X(\mu)\}$  is relatively compact in $\Hom
(W,T)/X (\mu)$ it suffices to show
that $\{\alpha_j+X(\mu_i)\}$  is relatively compact for all
$i$. In other words in proving the theorem we may assume
that $\theta (\mu)$ is a pure measure.

Let $E$ be the subset of $\widehat T$ consisting of all $\chi
\in \widehat T$ such that $\{\chi \circ \alpha_j\}$ is a bounded
subset of $\wh T$. Then $E$ is a subgroup of $\widehat T$. Also,  the
defining condition for $E$ readily implies that $\widehat T/E$ has
no nontrivial element of finite order. Let $S$ be the annihilator of $E$ in
$T$. Then $S$ is a compact subgroup, and the preceding observation
implies that $S$ is connected, namely a torus.
Let $S'$ be a subtorus of $T$ such that $T=S'S$, a direct product.
Then each $\alpha_j$ can be decomposed as $\beta_j+\gamma_j$
canonically, with $\beta_j\in \Hom (W,S')$, and $\gamma_j\in \Hom
(W,S)$, with $\Hom (W,S')$ and $\Hom (W,S)$ viewed as subgroups
of $\Hom (W,T)$ canonically.
%Correspondingly $\wt \alpha_j=\wt \beta_j \wt \gamma_j$ for
%all $i$. We shall show that $\{\wt \beta_j\}$ is relatively compact and
%$\wt \gamma_j\in A(\mu)$ for all $i$. This would prove the theorem.

%To prove that  $\{\wt \beta_j\}$ is relatively compact it suffices to
We note that $\{ \beta_j\}$ is relatively compact in $\Hom (W,S')$. To
see this it
suffices to know that $\{\chi \circ \beta_j\}$ is bounded  for all
characters $\chi$ on $S'$. Since the latter, viewed as characters on $T$,
belong to $E$ the desired property is immediate from the definition of $E$.

We now consider a Borel section $s$ of $L/T$ into $L$ such that
$\Lambda:= s(L/T)$
is locally compact, and identify $L$ and $\Lambda
\times T$ as Borel spaces. For $\varphi \in C_{c}(L/T)$ and $\chi \in
\widehat{T}$, we denote  by
$\varphi \otimes \chi$ the function on $L$ defined by $\varphi \otimes
\chi (x t)=\varphi (x) \chi (t)$, for all $x \in \Lambda$ and $ t \in
T$. The  set
of such functions is separating in $P(L)$, as is easily verified using
disintegration of measures on $L$ with respect to the projections of
$L$ on $L/T$. Hence, we can  test weak convergence of measures on $L$,
using these functions. We disintegrate $\mu$ as
$\mu=\int \delta_{x}\otimes \mu_{x} d\bar{\mu}(x)$ where $\mu_{x}\in
P(T)$ for all $x\in \Lambda$ and $\bar{\mu} \in P(L/T)$ is the
projection of $\mu$ on $L/T$.

Viewing elements of $\Hom (W,T)$ also as functions on $L/T$ via the
canonical projection of $L/T$ on to $L/R$,
for the measures $\tilde \alpha_j(\mu)$ we see that
$$\widetilde{\alpha}_{j}(\mu) (\varphi \otimes \chi)=\int \varphi (x)
\chi (t)  \chi(\alpha_{j} (x)) d\mu_{x} (t) d\bar{\mu} (x)=\int
\varphi (x) \widehat{\mu}_{x} (\chi) (\chi \circ \alpha_{j}) (x)
d\bar{\mu}(x)$$
for all $j$. We denote by $\mu^{\chi}_{\varphi}$ the projection on
$W=L/R$ of the
measure with density $\varphi(x) \widehat{\mu}_{x}(\chi)$ with
respect to $\bar{\mu}$. Then the value of the above integral is
$\widehat{\mu}^{\chi}_{\varphi} (\chi \circ \alpha_{j})$. Since
by the condition in the hypothesis a power $\theta (\mu)^d$ of $\theta
(\mu)$ has a density with respect to the Lebesgue measure on
$W$, the same holds for the measure
$\mu^{\chi}_{\varphi}$, and hence $\widehat{\mu}^{\chi}_{\varphi}$
vanishes at infinity on $\wh W$. Then, for $\chi \notin E $, and in
particular for any nontrivial character $\chi$ from $\wh S \subset \wh T$, $\{\chi \circ
\alpha_{j}\}$ is unbounded  on $\widehat W$, and hence the preceding
conclusion together with the fact that $\tilde \alpha_j(\mu) \to \tau
(\mu)$ implies that  $\tau
(\mu)(\varphi \otimes \chi)=0$. This means that the measure $\tau
(\mu)$ is translation-invariant under the action of elements of
$S$. Since $\tau$ is a continuous automorphism that pointwise fixes
$T$,  and in particular
$S$, it follows that $\mu$ is translation-invariant under the
$S$-action; thus $\mu_x$ is $S$-invariant for almost all $x\in L/T$. Since
for all $j$ the image of $\gamma_j$ is contained in $S$, this shows that
$\mu$ is invariant under $\wt \gamma_j$ for all $j$. Hence
$\gamma_j\in X (\mu)$ for all $j$. Since
$\alpha_j=\beta_j\gamma_j$ for all $j$ and $\{\beta_j\}$ is relatively
compact this shows that $\{\alpha_j+X(\mu)\}$ is relatively
compact.~\hfill $\Box$

\section{Some properties of roots of measures}

Given a Lie group $G$, a closed subgroup $H$ of $G$ and a subset
$\Lambda$ of $P(G)$ we say that
%{\it relatively compact modulo a closed subgroup} $H$ of $G$, and we write
$\Lambda/H$ {\it is relatively compact}, if for every sequence
$\{\lambda_i\}$ in $\Lambda$ there exists a sequence $\{h_i\}$  in $H$ such that
$\{\lambda_i h_i\}$ is relatively compact in $P(G)$.

\begin{theorem}\label{thm:cov}
Let $G$ be a Lie group (not necessarily connected) and suppose that
there exists a
continuous homomorphism $r:G \to \wt G$ onto an almost algebraic
group $\wt G$, such that $\ker r$ is  a discrete
central subgroup of
$G$. Let $\sigma \in P(G)$ and let $\Psi$ be a set of roots of
$\sigma$ in
$P(G)$. Let $M$ be the smallest closed normal $G(\Psi)$-invariant
subgroup of $Z^0(\sigma)$
such that $Z^0(\sigma)/M$ is a vector group and
the induced $G(\Psi)$-action on $ Z^0(\sigma)/M$ is trivial.
Then $\Psi/M$ is relatively compact.
\end{theorem}

\proof Let $\nu=r(\sigma)$. We first prove that $r(Z^0(\sigma))=Z^0(\nu)$: Let
$C=r^{-1}(Z(\nu))^0$. For any $g\in G(\sigma)$, $\{gzg^{-1}z^{-1}\mid
z\in C \}$ is a connected subset contained
in $\ker r$, and as the latter is discrete it consists only of the
identity element. Therefore $C$ is contained in $Z^0(\sigma)$. Also,
since $\ker r$ is discrete, $r(C)=Z^0(\nu)$. Thus we have
$Z^0(\nu)\subseteq r(Z^0(\sigma))$. The other way inclusion is easy to
see.

We next show that $\Psi /Z(\sigma)$ is relatively compact. Let
$\{\lambda_i\}$ be a sequence in $\Psi$. Then $\lambda_i$ are factors
of $\sigma$ and by a standard argument involved in the factor compactness
theorems (see \cite{DM-jtp}, \cite{DM-mz}, \cite{D-cimpa}) there exists a sequence
$\{x_i\}$ in $N(\sigma)$ such that $\{\lambda_ix_i\}$ and
$\{x_i^{-1}\sigma x_i\}$ are relatively compact in $P(G)$. Then
$\{r(x_i)^{-1}\nu r(x_i)\}$ is relatively compact in $P(\wt G)$. Since
$\wt G$ is almost algebraic this implies that $\{r(x_i)Z^0(\nu)\}$ is
relatively compact in $\wt N(\nu)/Z^0(\nu)$ (see~\cite{DM-mz}). Since
$r(Z^0(\sigma))=Z^0(\nu)$ this yields that $\{x_iZ^0(\sigma)(\ker r)\}$ is
relatively compact in  $N(\sigma)/Z^0(\sigma)(\ker r)$, and as  $Z^0(\sigma)(\ker
r)$ is contained in $Z(\sigma)$ it follows that  $\{x_iZ(\sigma)\}$ is
relatively compact in  $N(\sigma)/Z(\sigma)$. As $\{\lambda_ix_i\}$ is
relatively compact this implies in turn that there exists a
sequence $\{z_i\}$ in $Z(\sigma)$
such that $\{\lambda_iz_i\}$ is relatively compact. Therefore
$\Psi/Z(\sigma)$ is relatively compact.

Now let $\Phi=r(\Psi)$. We note that $\wt
G(\Phi)$ normalises $Z^0(\nu)$ and as they are almost algebraic
subgroups the product $\wt G(\Phi)Z^0(\nu)$ is also almost algebraic
and in particular a closed subgroup of $\wt G$.
%If $q: \wt G\to \wt G/\wt G_1$ is
%the quotient homomorphism then $q(\Phi)$ is a set of roots of $q(\nu)$
%in  $\wt G/\wt G_1$. As the latter is nilpotent $q(\Phi)$ is
%relatively compact. Hence we get that $\Phi/\wt G_1$ is relatively
%compact, and therefore to prove the theorem it suffices to show that
%$\Phi/Z(\nu, \Phi)$ is relatively compact.
%Since $Z^0(\nu)/M$ is compact
%this implies also that $\Phi/M$ is relatively compact.
Let $H=r^{-1}(\wt G(\Phi)Z^0(\nu))$. Let $\pi: H\to H/M$ be the quotient
homomorphism. Now, $\pi (Z(\sigma))$ is a closed
subgroup and as $\Psi/Z(\sigma)$ is relatively compact we get
that $\pi (\Psi)/\pi(Z(\sigma))$ is
relatively compact.

Let $Z'(\sigma)=Z^0(\sigma)(\ker r)$. Then $Z'(\sigma)$ is a closed subgroup
and $r(Z'(\sigma))=Z^0(\nu)$. Since
$Z^0(\nu)$ is of finite index in $Z(\nu)$ this implies that
$Z(\sigma)/Z'(\sigma)$ is finite. Hence we get that  $\pi (\Psi)/\pi(Z'(\sigma))$ is
relatively compact. From the definition of $M$ it follows that the
conjugation action of $G(\Psi)$ on $\pi
(Z^0(\sigma))$  is trivial. The subgroup of $\wt G$ consisting of
elements whose conjugation action on $\pi (Z^0(\sigma))$ is trivial is an
almost algebraic subgroup of $\wt G$. Since $\ker r$ is central and
$\wt G(\Phi)$ is the smallest almost
algebraic group containing $r(G(\Psi))$ it follows that the conjugation
action of $r^{-1} (\wt G(\Phi))$ on $\pi (Z^0(\sigma))$ is trivial. Since
$Z^0(\sigma)/M$ is abelian this shows that
$\pi (Z^0(\sigma))$ is contained in the center of $ H/M$, and hence so is
$\pi (Z'(\sigma))$. Together with
the fact that $\pi (\Psi)/\pi(Z'(\sigma))$ is
relatively compact, by Proposition~3.4 of \cite{DM-inv}
we get  that $\{\pi (\Psi)\}$ is relatively compact. Since
$\Psi/Z(\sigma)$ is relatively compact and $\ker \pi
=M$ this shows that $\Psi /M$ is relatively
compact; for a sequence $\{\lambda_i\}$ in $\Psi$ if $\{z_i\}$ is a
sequence in $Z(\sigma)$ such that $\{\lambda_iz_i\}$ is relatively
compact then we get that $\{\pi (z_i)\}$ is relatively compact and so
$\{z_i\}$ can be replaced by a sequence in $\ker \pi$. This proves the
theorem.
\hfill $\Box$

\medskip
As before let $Q$ be the smallest closed normal subgroup of $\wt G(\nu)$ such
that $\wt G(\nu)/Q$ is a vector group.
 We note that $Q$ is the subgroup of $\wt
G(\nu)$ containing $[\wt G(\nu), \wt G(\nu)]$, such that $Q/[\wt
G(\nu), \wt G(\nu)]$ is the maximal compact subgroup of $\wt G(\nu)/[\wt
G(\nu), \wt G(\nu)]$. Hence $Q$ is an almost algebraic
subgroup of $\wt G$. Moreover, since $\wt G(\nu)/Q$ is a vector group, and
in particular has no finite subgroup, it follows that $Q$ is
Zariski closed in $\wt G(\nu)$. Therefore $\wt G(\nu)/Q$ is
an almost algebraic group. Let $\wt V= \wt G(\nu)/Q$, equipped with its
structure as an algebraic group. The characteristic property of $Q$ as
a subgroup of $\wt G(\nu)$ shows also that $Q$ is normal in $\wt
N(\nu)$. Let $\eta:\wt N(\nu)
\to \wt N(\nu)/Q $ be the canonical quotient homomorphism; we have
$\eta (\wt G(\nu))=\wt V$. For $g\in \wt
N(\nu)$ let $c_g:\wt V\to \wt V$ be the transformation
induced by the conjugation action of $g$ on $\wt G(\nu)$. Let $\wt
N_1(\nu)$ be the subgroup of $\wt
N(\nu)$ consisting of all elements $g$ such that the measure $\eta
(\nu)$  is invariant under the action of $c_g$.
%Let $I=\{c_g\mid g\in \wt N_1^0(\nu)\}$.
% Also, for any set of roots $\Phi$ of $\nu$ let $I(\Phi)=\{c_g
%\mid g\in \wt G(\Phi)\}$.

\begin{remark}\label{rem:inv}
{\rm
Given a  Lie group $G$ (or more
generally a locally compact group) and $\mu \in P(G)$, if
$G(\mu)$ is abelian and $\lambda$ is a root of $\mu$
  supported on $xG(\mu)$, with $x\in N(\mu)$, then $\mu$ is
  invariant under the conjugation action of $x$ on $G(\mu)$ (see
 \cite{DS}, Lemma~2.2).
}
\end{remark}

\begin{proposition}\label{cpt}
Let the notation be as above. Then we have the following;

i) $\{c_g\mid g\in \wt N_1^0(\nu)\}$ is a compact subgroup of
$GL(\wt V)$;

ii) if $\Phi$ is a set of roots of $\nu$ in $P(\wt
G)$ then $\{c_g \mid g\in \wt G(\Phi)\}$ is compact and preserves $\eta(\nu)$.
\end{proposition}

\proof i) This follows from Corollary~2.5 of \cite{D-isr},
since $\{c_g\mid g\in \wt N(\nu)\}$ is an almost algebraic subgroup of
$GL(V)$ and $\wt G(\eta (\nu))=\wt V$.

\iffalse
As $\wt V$ is an almost
algebraic group it can be
expressed as $\wt V_s\oplus \wt V_u$ where $\wt V_s$ and
$\wt V_u$ are almost algebraic subgroups consisting of semisimple elements
and unipotent elements respectively.  Then for all  $g\in \wt
N_1(\nu)$, $\wt V_s$ and $\wt V_u$ are
invariant under $c_g$, and furthermore, since
the group of algebraic automorphisms of $\wt V_s$ is countable it follows
that the action of $c_g$ on $\wt V_s$ is trivial for $g\in \wt
N_1^0(\nu)$. The restrictions of
$c_g$,  $g\in \wt N_1^0(\nu)$, to $\wt V_u$ form an almost algebraic and
hence closed subgroup of $GL(\wt V_u)$, leaving invariant the image of $\eta
(\nu)$ on $\wt V_u$. Since $\wt V_u$ consists of unipotent elements every
vector subspace of $\wt V_u$ is an almost algebraic subgroup. Hence
the image of $\eta (\nu)$  on $\wt V_u$ is a measure not supported on any
proper vector subspace. It is well-known that this condition implies
that the elements of $GL(\wt V_u)$ leaving it
invariant form a compact subgroup; the reader may refer to Corollary~2.2
in \cite{D-cimpa}, if necessary.
Thus in particular the restrictions of
$c_g$,  $g\in \wt N_1^0(\nu)$, to $\wt V_u$ form a compact
subgroup. Hence
we get that  $\{c_g\mid g\in \wt N_1^0(\nu)\}$ is a compact subgroup.
\fi

ii) For all $\rho \in \Phi$, $\eta (\rho)$ is a root of $\eta (\nu)\in
P(\wt V)$ and as $\wt V$ is abelian, by Remark~\ref{rem:inv}
 $\supp \rho$ is contained in  $\wt N_1(\nu)$. The compactness
of $\{c_g\mid g\in \wt N_1^0(\nu)\}$ implies in particular that
$\wt N_1(\nu)$ is an almost algebraic
subgroup. Therefore  $\wt G(\Phi)$  is contained in $\wt N_1(\nu)$. Also,
since $\wt G(\Phi)$ is an almost algebraic subgroup $\{c_g \mid
g\in \wt G(\Phi)\}$ is closed. Hence by (i) it is compact, and it
preserves $\eta (\nu)$, since  $\wt G(\Phi)$  is contained in $\wt
N_1(\nu)$.
%The second
%assertion readily follows from the continuity of $g \mapsto c_g$ and
%fact that $c_g$ is trivial for all $g\in \wt G(\nu)$.
\hfill $\Box$

\section{Embedding of infinitely divisible measures}

Our next objective is to prove the following theorem, from which
Theorem~\ref{thm:main} can be deduced.

For any finite measure
$\lambda$ on a Lie group $G$ we denote by $I(\lambda)$ the
subgroup $\{g\in G\mid g\lambda = \lambda g\}$ and by
$J(\lambda)$ the subgroup $\{g\in G\mid g\lambda
=\lambda
g=\lambda\}.$ We note that $I(\lambda)$ and $J(\lambda)$ are closed
subgroups and $J(\lambda)$ is normal in
$I(\lambda)$.

\begin{theorem}\label{thm:main-tech}
Let $G$, $T$ and $q:G\to G/T$ be as in Theorem~\ref{thm:main}. Let
$\mu\in P(G)$,  $\sigma=q(\mu)$, and suppose that $\sigma$ has no
nontrivial idempotent
factor in $P(G/T)$. Let $\Psi$ be a root cluster of $\mu$,
 and let $M$ be the smallest closed normal  $ G(\Psi)$-invariant
subgroup of $Z^0(\sigma)$
such that $Z^0(\sigma)/M$ is a vector group and
the induced $G(\Psi)$-action on $ Z^0(\sigma)/M$ is trivial.
Suppose that,  in $G/T$, $G(q(\Psi))$ normalizes  $q(I(\mu))\cap M$.
Then there exists  a  rational embedding $\{\mu_r\}_{r\in \Q^+}$ of
 $\mu$ such that $\mu_r\in \Psi$ for all $r\in \Q^+$.
\end{theorem}

 Apart from the results of the previous
sections, for the proof of Theorem~\ref{thm:main-tech}  we need certain
general results. We begin by recalling these.

Firstly we recall the following result from \cite{D}, extending
analogous results proved in \cite{DM-adv} under some restrictions on the
subgroup $N$ as in the statement of the theorem below. Let $H$ be a
Lie group and $\sigma \in P(H)$.
We recall (from \cite{D}) that a root $\rho$ of
$\sigma$  is said to be {\it compatible}  with a Lie subgroup $N$ of
$N(\sigma)$ if $\rho\in P(xG(\sigma))$ for some $x\in N(\sigma) $ which
normalises $N$. The following is Corollary~2.4 from \cite{D}, with a
co-finite subset of  $\N$ in the place of $\Sigma$ there.

\begin{theorem}\label{thm:roots}
Let $H$ be a Lie group.
Let $\sigma \in P(H)$ and for each $k$ let $\rho_k$ be a $k!\,$ th
root of $\sigma$. For all $n\in \N$ let $S_n=\{\rho_k^{k!/n!}\mid k\geq n\}$.
Let $N$ be a closed compactly generated  nilpotent subgroup of $Z(\sigma)$.
Suppose that all $\rho_k$ are compatible with $N$ and that $S_n/N$ is
relatively compact for all large $n$.
Then there exist sequences $\{k_j\}$ in $\N$,  $\{z_j\}$ in $N^0$
and $n_0\in \N$ such that $k_j\to \infty$ as $j\to \infty$ and
$\{z_j\rho_{k_j}^{k_j!/n!}z_j^{-1}\mid j\in \N, k_j! \geq n!\}$ is relatively
compact for all $n\geq n_0$.
\end{theorem}

We need also the following result from \cite{S}. It may be noted
that the condition in the hypothesis of Theorem~\ref{thm:main-tech} is
connected with the use of this theorem.

\begin{theorem}\label{thm:riddhi}
Let $G$ be a Lie group and $\sigma \in P(G)$.
Let $\lambda_j\in P(G)$ be a $n_j$th root of $\sigma$, where $n_j\to
\infty$ as $j\to \infty$. Let $\{x_j\}$ be a sequence in $Z(\sigma)$
such that $\{\lambda_jx_j\}$ converges to a measure $\lambda \in
P(G)$.  Then there exists $x \in I(\sigma)$ such that
the support of $\lambda $ is contained in $xJ(\sigma)$.
\end{theorem}

\proof Each $\lambda_jx_j$ is a factor of $\sigma$ and hence it follows
that the limit $\lambda$ is a factor of $\sigma$. Also, for any $k\in
\N$, $(\lambda_jx_j)^k$ can be expressed as $\lambda_j^ky_j^{(k)}$,
with $y_j^{(k)}\in Z(\sigma)$ (see \cite{DM-adv}, Proposition 2.2 for an
idea of the proof). Hence $(\lambda_jx_j)^k$ is a
factor of $\sigma$ for all $j$ such that $n_j\geq k$. As
$(\lambda_jx_j)^k\to \lambda^k$ as $j\to \infty$ and $n_j\to \infty$,
we get that $\lambda^k$ is a factor of $\sigma$ for all $k$. The
theorem now follows from Theorem~2.4 of~\cite{S}.~\hfill $\Box$

\medskip
To prove the general case of Theorem~\ref{thm:main-tech} we will need
also the following theorem, which enables to reduce to a situation
where $Z^0(\sigma)$
contains a cocompact simply connected nilpotent normal subgroup.
We postpone proving it, until \S~7, and go over to the
proof of Theorem~\ref{thm:main-tech}; we note that for certain
Lie groups, such as the groups involved in Corollary~\ref{walnut}, the
conclusion of Theorem~\ref{reduction} is readily seen to hold.

\begin{theorem}\label{reduction}
Let $\wt G$ be an almost algebraic group and $\nu \in P(\wt G)$. Let
$\Phi$ be a set of roots of $\nu$. Suppose that $\wt G (\Phi)/\wt G(\nu)$ is
connected. Let $K$ be the
smallest almost algebraic normal subgroup of $Z^0(\nu)$ such that
$Z^0(\nu)/K$ is a vector group.
Then there exists an almost algebraic
subgroup $\wt H$ of $\wt G(\Phi)$ such that the following conditions
are satisfied:

i) for all $\rho \in \Phi$ there exists $z\in K$ such that $z\rho
z^{-1}\in P(\wt H)$, and

ii)  if $U$ is  the unipotent radical of $\wt
G^0(\Phi)$ and $P=\wt G(\nu)U$, then $\wt H$ contains $P$ and
the quotient $(Z^0(\nu)\cap \wt H)/(Z^0(\nu)\cap P) $ is compact;
furthermore $Z^0(\nu)\cap \wt H$ contains a simply connected nilpotent
cocompact normal almost algebraic subgroup.

\end{theorem}

We now begin the proof of Theorem~\ref{thm:main-tech}. Let the
notation be as in the hypothesis. Also, let $\nu =p(\mu)$. Let
$\Psi$ be a root cluster of $\mu$ in $P(G)$.
The subgroup $\wt G(p(\Psi))$ (see \S~2 for definition) being an almost
algebraic subgroup has only finitely many connected components; let
$c$ be the number of the components.
Now, let $\Psi'=\{\rho^c\mid \rho, \rho^c  \in \Psi\}$. It is straightforward
to see that $\Psi'$ is also a root cluster of $\mu$ contained in
$\Psi$. Furthermore  $\wt G(p(\Psi'))/\wt G(\nu)$ is
connected. Therefore, in proving
Theorem~\ref{thm:main-tech}, replacing $\Psi$ by $\Psi'$ and modifying the
notation we may assume that  $\wt G(p(\Psi))/\wt G(\nu)$ is connected.

The conditions of Theorem~\ref{reduction} are satisfied for $\wt G$,
$\nu$  as above and $\Phi=p(\Psi)$. Let $\wt H$ be an almost algebraic
subgroup
of $\wt G(\Phi)$ for which the assertions as in that theorem are
satisfied.
Now let $G'=p^{-1}(\wt H)/T$, and let $q:p^{-1}(\wt H)\to G'$ and
$r:G'\to \wt H$ be the quotient
homomorphisms. We note that $r$ is a covering homomorphism onto the
almost algebraic group $\wt H$. Now let $\Psi^*=\Psi \cap P(p^{-1}(\wt
H))$. We note that $\Psi^*$ is a root cluster for $\mu$, which we now
view as a measure on $p^{-1}(\wt H)$: given $n\in \N$ there exists
$\rho \in \Psi$ such that $\rho^n=\mu$ and by
condition~(i) in the conclusion of Theorem~\ref{reduction} there
exists $x\in p^{-1}(K)$ such that $x\rho x^{-1}\in P(p^{-1}(\wt H))$;
since $p^{-1}(K)\subset Z(\mu)$ and $\Psi$ is a root cluster we get
that $x\rho x^{-1}\in \Psi$, showing that $\Psi^*$ contains a $n$ th
root of $\mu$, for any $n$; the other conditions as in the definition
of a root cluster
are straightforward to verify. Replacing $\Psi$ by $\Psi^*$ if
necessary and modifying the notation, in proving
Theorem~\ref{thm:main-tech} we may without loss of generality assume that
$\Psi$ is contained in $p^{-1}(\wt H)$.

Now let $M$ be the subgroup as in the hypothesis of
Theorem~\ref{thm:main-tech}.
Then by Theorem~\ref{thm:cov} $q(\Psi)/M$ is
relatively compact. By assertion~(ii) of
Theorem~\ref{reduction} $Z^0(\nu)\cap \wt H$ admits a
cocompact almost algebraic
simply connected nilpotent normal subgroup. We note that $r(M)$ is an
almost algebraic subgroup of $Z^0(\nu)\cap \wt H$ and hence we get
that $r(M)$ contains  a cocompact almost algebraic
simply connected nilpotent normal subgroup.
Since such a subgroup is
necessarily unique, it follows that it is normalised by $\wt G(\Phi)$.
Since $r(Z^0(q(\mu)))=Z^0(\nu)\cap \wt H$ this implies that $M$
contains  a cocompact
simply connected nilpotent normal subgroup, invariant under the action
of $\wt G(\Phi)$; we shall denote the subgroup, which is unique, by
$N$. Then $q(\Psi)/N$ is relatively compact.
%Let $N'=\{gT\in G'\mid g\in p^{-1}(N), g\mu g^{-1}=\mu\}$.

\bigskip
Towards the proof of the Theorem~\ref{thm:main-tech} we shall now first prove
the following.

\begin{proposition}\label{thm:key}
Let the notation be as above. Also, for all $k\in \N$ let
$\rho_k$ be a $k!$~th root of $\mu$ contained in $\Psi$.
For all $n\in \N$ let  $S_n=\{q(\rho_k)^{k!/n}\mid k\geq
n\}\subseteq P(G')$. Let $N'=N\cap q(I(\mu))$  Then $S_n/N'$ is
relatively compact for all large $n$.

\end{proposition}

\proof
%Before going over to the proof of the theorem we shall describe a
%subgroup of $N'$ needed in the proof.
Let $R$ be the subgroup  of
$p^{-1}(\wt G(\nu))$ consisting of all $g$ such that $\psi_{p(z)}(p(g))=e$,
the identity element, for all $z\in q^{-1}(M)$. By Lemma~\ref{psi} $R$ contains
$p^{-1}(Q)$. Therefore $p^{-1}(\wt G(\nu))/R$ is a quotient of $\wt
G(\nu)/Q$ and hence a vector group. Let $W=p^{-1}(\wt
G(\nu))/R$, equipped  with the induced vector space structure.
Also, $p(R)/Q$ is invariant under the action of $\wt G(\Phi)$ on $\wt
G(\nu)/Q$. Moreover by Corollary~\ref{cor:nofp} every fixed point of
the $\wt G(\Phi)$-action on $\wt G(\nu)/Q$ is contained in $p(R)/Q$.
Since the $\wt G(\Phi)$-action on $\wt G(\nu)/Q$ is through a compact group
(see Proposition~\ref{cpt}) this further implies that the
induced $\wt G(\Phi)$-action on $W$ has no nonzero fixed point.

Let $\theta: p^{-1}(\wt G(\nu))\to W$ be the quotient homomorphism. We
write $\mu$ as $\mu_0+\mu'$, where $\mu_0$ is
supported on $R$ and $\mu' (R)=0$.

Now let $\Delta=\{(k,n)\mid k,n\in \N, k\geq n\}$. For any
$\delta=(k,n)\in \Delta$ let $\theta_\delta =q(\rho_k)^{k!/n}$.
As $q(\Psi)/N$ is relatively compact we get that there
exists a family $\{x_\delta\}_{\delta \in \Delta}$ in $ N$ such
that $\{\theta_\delta
x_\delta\}$ is relatively compact. Moreover for a subset $E$ of
$\Delta$ and a closed normal subgroup $H$ of $ N$, $\{\theta_\delta\mid
\delta \in E\}/H$ is relatively compact if and only if
$\{x_\delta H\}_{\delta \in E}$ is relatively compact in $N/H$.
We note that  $\psi_z$ is trivial for all $z\in [N,N]$, and hence
it follows that $[N,N]$ is contained in $q(I(\mu))$. Therefore
$N\cap q(I(\mu))$ is a normal subgroup of $N$.
To prove the proposition it therefore suffices to show that
$\{x_{(k,n)} \mid k\geq n \}/(N\cap q(I(\mu)))$ is relatively
compact  for all large $n$. By the definition
of $R$, $\psi_{p(z)}(p(g))=e$ for all $g\in R$ and $z\in q^{-1}(M)$, and in
particular for all $z\in q^{-1}(N)$.  Since the support of $\mu_0$ is
contained in $R$, this implies that $q(I(\mu_0))$ contains $N$. Hence
$N\cap q(I(\mu))=N\cap q(I(\mu'))=N'$, say.

Suppose, if possible that $\{x_{(k,n)}N'\mid k\geq n\}$ is
not relatively compact for infinitely many $n$. We can then
find sequences $\{k_j\}$ and   $\{n_j\}$ in $\N$, such that, $k_j\geq
n_j$, $n_j\to \infty$ and
$\{x_{(k_j,n_j)}N'\}$ has no convergent subsequence in $
N/N'$. For all $j\in \N$ let $\theta_j=q(\rho_{k_j}^{k_j!/n_j})$ and
$x_j=x_{(k_j,n_j)}$. Then $\{\theta_j x_j\}$ is
relatively compact. Let $\lambda \in P(G')$ be a limit point of
$\{\theta_jx_j\}$. By Theorem~\ref{thm:riddhi}
there exists  $x \in I(q(\mu))$ such that  $\supp \lambda$ is
contained in $xJ(q(\mu))$ (notation as in Theorem~\ref{thm:riddhi}). The
 normalised Haar measure of $J(q(\mu))$ is an idempotent factor of
$q(\mu)$, and  since by hypothesis $q(\mu)$ has no nontrivial idempotent
factor  we see that
$J(q(\mu))$ is trivial. Hence  $\lambda =\delta_x$ with
$x \in I(q(\mu))$.

For every $j$ let $\phi_j$ be a $n_j$ th root of $\mu$ such that
$q(\phi_j)=\theta_j$. Let  $y_j \in G$ be such
that $q(y_j)=x_j$. Then $\{\phi_jy_j\}$ has a limit point, say
$ \alpha \in P(G)$ such that $q( \alpha)=\lambda=\delta_x$.
Let $y\in G$ be such that $q(y)=x$. Then
 $ \alpha =y \beta$ for some $ \beta \in P(T)$.

For all $j$  we have the identity $\mu
=(\phi_jy_j)(y_j^{-1}\phi_j^{k_j-1})$. Since
$\{\phi_jy_j\}$ is relatively compact it follows also
that $\{y_j^{-1}\phi_j^{k_j-1}\}$ is relatively compact
(cf. \cite{P}, Chapter~III, Theorem~2.1) and hence has a
subsequence which converges, to say $\gamma \in P(G)$. Then passing to
limit along such a subsequence we have $\mu
=y  \beta \gamma$. Also
$\{y_j^{-1}\mu y_j\}=\{(y_j^{-1}\phi_j^{k_j-1})(\phi_jy_j)\}$
converges along the subsequence to
$ \gamma (y \beta)=\gamma \beta y =\beta \gamma y= y^{-1} \mu y$; note
that since
 $\beta \in P(T)$, $y\beta =\beta y$ and $\gamma \beta =\beta
 \gamma$. Thus $y^{-1}\mu y$ is a limit of a subsequence of
$\{y_j^{-1}\mu y_j\}$.  It follows also that $\{y_j^{-1}\mu' y_j\}$
has a subsequence converging to $y^{-1} \mu' y $.
%Recall that $\widehat \lambda_i$ vanishes at infinity on $\widehat W_i$.

Now let  $L=p^{-1}(\wt G(\nu))$. We  apply Theorem~\ref{thm:FT} to
this $L$, with $W$ as above, and $\mu'$ in the place of $\mu$ there.
 Let $\lambda =\theta (\mu')$.  We recall that the $\wt
 G(\Phi)$-action on $W$ has no nonzero fixed point, and it preserves
 $\lambda$ (see Proposition~\ref{cpt}).  Then  by Theorem~\ref{yves},
the pure components $\lambda_i$ of $\lambda$ are such that
$\lambda_i^{d_i}$ have a density on the subspace spanned by $\supp \lambda_i$,
for some $d_i$, and
%$\widehat \lambda $ vanishes at infinity, and
hence the
condition in the hypothesis of   Theorem~\ref{thm:FT} is satisfied.
We note also for any $z \in q^{-1}(N)$ and $g\in L$ we have
$z^{-1}gz =g \alpha_z(\theta (g))$, where $\alpha_z:W\to T$ is the map
given by $\alpha_z(\theta (g))=g^{-1}z^{-1}gz=\psi_{z^{-1}}(p(g))$ (in
  the notation as before) for
  all $g\in L$. Since  $\{y_j^{-1}\mu' y_j\}$
has a subsequence converging to $y^{-1} \mu' y $, in view of these
observations Theorem~\ref{thm:FT} implies that the sequence
$\{x_jN'\}$ has a convergent
subsequence in $N/N'$. As this contradicts the choice of
$\{x_j\}$, we get that
 $S_n/N'$ is relatively compact for all large $n \in \N$. This
 proves the proposition.~\hfill $\Box$

\bigskip
%\pagebreak
\noindent{\it Completion of the proof of Theorem~\ref{thm:main-tech}}:
We follow the notation as above. By Proposition~\ref{thm:key} $S_n/N'$
is relatively compact for all large $n$. Recall that by hypothesis
$q(I(\mu))\cap M$ is normalized by $\wt G(\Phi)$. Also, as seen above
$N$ is normalized by $\wt
G(\Phi)$. Hence $N'=N\cap (q(I(\mu))\cap M)$ is normalized by
$\wt G(\Phi)$. This shows that all
$q(\rho_k)$, $k\in \N$ are compatible with $N'$.  Hence by
Theorems~\ref{thm:roots} there exist
sequences $\{k_j\}$ in $\N$,  $\{z_j\}$ in $N'$
and $n_0\in \N$ such that $k_j\to \infty$ as $j\to \infty$, and
$\{z_jq(\rho_{k_j})^{k_j!/n!}z_j^{-1}\mid j\in \N, k_j! \geq n!\}$ is relatively
compact for all $n\geq n_0$.
%sequence $\{z_k\}$ in $N$ such that $\{z_k\rho_k^{k!/n!}z_k^{-1}\}$ is
%relatively compact for all large $n$, say $n\geq n_0$.
For all $j$ let $\zeta_j \in q^{-1}(z_j)$.
Since the conjugation action of $q^{-1}(N')$ leaves $\mu$ invariant,
each $\zeta_j\rho_{k_j}^{k_j!/n!}\zeta_j^{-1}$ is a root of $\mu$.
For all $n\geq n_0$, since
$\{q(\zeta_j\rho_{k_j}^{k_j!/n!}\zeta_j^{-1})\}=\{z_jq(\rho_{k_j})^{k_j!/n!}z_j^{-1}
\}$ is
relatively compact and the kernel of $q$ is compact
it follows that $\{\zeta_j\rho_{k_j}^{k_j!/n!}\zeta_j^{-1}\}$ is
relatively compact.  Following a standard argument with limits we can
now find $k!$ th roots $\mu_{1/k!}$, $k\geq 2$, such that
$(\mu_{1/k!})^{k}=\mu_{1/(k-1)!}$ for all $k$. We now define
for $r=p/q$, where $p,q\in \N$, $\mu_r$ to be
$\mu_{1/q!}^{p(q-1)!}$; in view of the preceding conditions it follows
that $\mu_r$ defined in this way is indeed independent of the
representation of $r$ as $p/q$. It is also straightforward to see that
$\{\mu_r\}$ is a rational embedding of $\mu$ contained in $\Psi$.
 This proves the first statement in the theorem. As noted before the
 later statements are well-known consequences of the first one. This
completes the proof of the theorem (assuming the validity of
Theorem~\ref{reduction}).~\hfill $\Box$

\section{Measures on almost algebraic groups}

We now prove Theorem~\ref{reduction} on ``reducing'' an almost
algebraic group so that $Z(\nu)$  contains a simply
connected nilpotent cocompact normal subgroup. The proof is along the
lines of the reduction theorem in \cite{DM-inv} and \cite{DM-adv},
the crucial additional point being that the conjugating elements are
chosen from the subgroup $K$, rather than the whole of $Z(\nu)$.

\medskip
\noindent {\it Proof of  Theorem~\ref{reduction}}:  Let $U$ and $P$ be
as in (ii)~in the statement of the theorem. As $\wt G(\nu)$ and $U$
are normal almost algebraic subgroups of $\wt G(\Phi)$, so is $P$.
As $U$ is the unipotent radical of $\wt G^0(\Phi)$, $\wt
G(\Phi)/P$ is a reductive almost algebraic group. Also, since  $\wt G
(\Phi)/\wt G(\nu)$ is connected, so is $\wt
G(\Phi)/P$. Now let  $Z'=(\wt G(\Phi) \cap Z(\nu))^0P$, which is an
almost algebraic subgroup of $\wt G(\Phi)$.  Then $Z'/P$ is a
connected almost algebraic normal
subgroup of $\wt G(\Phi)/P$ and hence there exists an almost algebraic
normal subgroup
$R$ of $\wt G(\Phi)$ containing $P$ and such that $\wt G(\Phi)/P$ is
the almost direct product of $R/P$ and $Z'/P$, namely,
$\wt G(\Phi)=RZ'$ and $(R\cap Z')/P$ is finite. Now let $\wt H$ be a
closed subgroup of $\wt G(\Phi)$ containing $R$ and such that $\wt H/R$ is
a maximal compact subgroup of $\wt G(\Phi)/R$. Let $K'=K\cap Z'$. Then
$K'$ is an almost algebraic normal subgroup of $Z'$ and $Z'/K'$ is a
vector group. We note that $RK'$ is an almost algebraic normal
subgroup of $\wt G(\Phi)$ such that  $\wt G(\Phi)/RK'$ is a vector
group, and hence in particular $\wt H$ is contained in $RK'$. Clearly
$\wt H$ is a maximal compact subgroup of $RK'$.

Let $\rho \in \Phi$ and $x\in \wt G(\Phi)$ be such that
$\rho \in P(x \wt G(\nu))$. As $x \wt G(\nu)$  is of
finite order in $\wt G(\Phi)/\wt G(\nu)$, so is $xR$ in
 $\wt G(\Phi)/R$. By conjugacy of maximal compact subgroup we
therefore get that $xR$ is conjugate in  $K'R/R$ to an element
 of $\wt H/R$. Therefore there exists a
 $z\in K'$ such that  $zxz^{-1}\in \wt H$.
 Hence $z\rho z^{-1}\in  P(\wt H)$. This proves~(i).

Now,  $(Z^0 (\nu))\cap \wt H)P$ is an almost algebraic and hence a
closed subgroup, and $(Z^0 (\nu))\cap \wt H)P/P$ is topologically
isomorphic to $(Z^0 (\nu)\cap \wt H)/(Z^0 (\nu)\cap P)$.
%$(\wt G (\nu)\cap Z^0(\nu))(U\cap Z^0(\nu))=(\wt
%G(\nu)U)\cap Z^0(\nu)=P\cap Z^0(\nu)$.
Therefore to prove~(ii) it suffices to show that $(\wt H\cap Z^0 (\nu))P/P$ is
compact. Recall that $\wt G(\Phi)/P$ is an almost direct product of
 $Z^0(\nu)P/P$ and $R/P$. Also $\wt H$ contains $R$
 and $\wt H/R$ is compact. The preceding observation therefore implies
 that $(\wt H\cap Z^0 (\nu))P/P$ is compact, as sought to be shown.

Now let $W$ be the unipotent radical of $Z(\nu)\cap P$. We show
that  $(Z(\nu)\cap P)^0$ is contained in
$(Z(\nu) \cap \wt G(\nu))W$. As $W$ is the unipotent radical
of $Z(\nu)\cap P$ by Levi decomposition there exists a
reductive almost algebraic subgroup $L$ such that $(Z(\nu)\cap
P)^0=LW$. Since $U$ is contained in the
unipotent radical of $P^0$, we get that there exists $u\in U$ such that
$uLu^{-1}$ is contained in $\wt G^0(\nu)$. But
$U$ is contained in $\wt N(\nu)$, so this means that $L$ is contained
in $\wt G(\nu)$.  This proves the claim as above.
Now let $A$ be the maximal almost algebraic
 vector subgroup of the abelian almost algebraic group $\wt G(\nu)\cap
 Z(\nu)$; in particular $\wt G(\nu)\cap
 Z(\nu)/A$ is compact. Then  $AW$ is a  normal
 almost algebraic subgroup of $Z(\nu)\cap \wt H$ and since $A$
 centralizes $W$,  $AW$ is a
simply connected nilpotent Lie group. Also, since $Z^0(\nu)\cap  P$
is cocompact in $Z^0(\nu)\cap \wt H$ the above observation shows that
$(Z^0(\nu)\cap \wt H)/AW$ is compact. This proves the theorem. \hfill $\Box$

\section{Concluding results}

In this section we shall deduce the results stated in the introduction.
To relate Theorem~\ref{thm:main} to Theorem~\ref{thm:main-tech} we
need the following proposition.

\begin{proposition}\label{prop:full}
Let the notation $G,\wt G,p,q$ and $T$ be as in
Theorem~\ref{thm:main}. Let $\mu\in P(G)$, $\sigma=q(\mu)
\in P(G/T)$ and $\nu= p(\mu) \in P(\wt
G)$.  Let $\Psi$ be a root
cluster of $\mu$ such that  $\wt G(p(\Psi))/\wt G(\nu)$ is connected.
Let $M$ be the subgroup of $Z^0(\sigma)$ as in
Theorem~\ref{thm:main-tech}.
Then $q(I(\mu)) \cap M$ is normalized by $G(q(\Psi))$.
\end{proposition}

\proof  Let $R$ be the subgroup of
$p^{-1}(\tilde G(\nu))$ consisting of all $g$ in $G$ such that
$\psi_{p(z)}(p(g))=e$ for all $z\in q^{-1}(M)$. From Remark~\ref{invrem} we see
that $R$ is a closed normal
subgroup of $p^{-1}(\tilde G(\nu))$ and $p^{-1}(\tilde G(\nu))/R$ is a
vector group, a quotient of $\wt V$, which we shall denote by
$W$. Also, $R$ is
invariant under the conjugation action of $G(\Psi)$ and by
Corollary~\ref{cor:nofp} the induced $G(\Psi)$-action on
$p^{-1}(\tilde G(\nu))/R$ has no nonzero fixed point. Let
$\Phi=p(\Psi)$ and $G_1=p^{-1}(\wt G(\Phi))$; we note that $G_1$
contains  $p^{-1}(\tilde G(\nu))$ and $R$ as normal subgroups.
Let $\theta:G_1\to G_1/R$ be the quotient homomorphism; then
$\theta(p^{-1}(\tilde G(\nu)))=W$.

We write $\mu$  as $\mu_0+\mu'$ where $\mu_0\in P(R)$ and
$\mu'(R)=0$. As $q^{-1}(M)$ centralizes $R$
it follows that $I(\mu)\cap q^{-1}(M)=I(\mu')\cap q^{-1}(M)$. Now we express $\mu'$
as $\mu'=\sum_ia_i \mu_i$, with $i$ running over a suitable countable
indexing set, $a_i$'s positive real numbers such that $\sum_ia_i=1$,
and $\mu_i$ are probability measures such that $a_i\theta (\mu_i)$
are the pure  components of $\theta (\mu)$; (see \S~3). Such a
decomposition is
unique and in particular this implies that $I(\mu)\cap q^{-1}(M) =\cap_i
(I(\mu_i)\cap q^{-1}(M))$.
For each $i$ let $W_i$ be the subspace of $W$ spanned by
$\theta(\mu_i)$, and let  $d_i=\dim W_i$. By Proposition~\ref{cpt},
$\theta (\mu)$ is invariant under the induced action of $\wt G(\Phi)$
on $W$. The action of $\wt G(\nu)$ on $W$ is trivial, and hence the
$\wt G(\Phi)$-action involved is through the quotient group $\wt
G(\Phi)/\wt G(\nu)$. As the latter is a connected group, by
Lemma~\ref{inv-comp} we get that $a_i\theta (\mu_i)$, which
are the pure components of $\theta (\mu)$, are $\wt G(\Phi)$-invariant
for all $i\in I$. In particular each $W_i$ is also $\wt
G(\Phi)$-invariant, and
in view of the invariance of the measure the action on $W_i$ is via a
compact subgroup of $GL(W_i)$.
Furthermore since the $\wt G(\Phi)$-action on $W$ has no
nonzero fixed point and by hypothesis $\wt G(\Phi)/\wt G(\nu)$ is
connected, by Theorem~\ref{yves} we get that
 $\theta(\mu_i)^{d_i}$ has a density in $W_i$.

Now consider the partition of $G$ into cosets of $T$ and let
$\{\mu_\xi\}_{\xi\in G/T}$ be a system of conditional measures
 for $\mu$ with respect to the partition. It is easy to see that there exists
a sequence of pairwise disjoint Borel subsets $\{B_i\}$ of $G/T$ such
that $q(\mu_i)(B_i)=1$. Therefore the system of conditional
measures as above can also be chosen so that it is a system of
conditional measures for each $\mu_i$, over the quotient measure
$q(\mu_i)$, for each index $i$.
For $g\in q^{-1}(M)$, each fibre
 $\xi= xT$, $x\in G$, is
 invariant under  the conjugation action of $g$
and the action on $\xi$ is given by translation by the element
 $g_\xi=gxg^{-1}x^{-1}\in T$, which is independent of the choice of
 $x$ in $\xi$. Hence for $g\in I(\mu)\cap q^{-1}(M)$, for $\sigma$-almost all
 $\xi$, $\mu_\xi$ is invariant under the translation
 action of $g_\xi$. For $\xi\in G/T$ let $S_\xi$ be the smallest
 closed subgroup of $T$ containing
 $\{g_\xi\mid g\in I(\mu)\cap q^{-1}(M)\}$. Then for $\sigma$-almost all
 $\xi$, $\mu_\xi$ is invariant under $S_\xi$.
We shall show that $S_\xi$ is constant $q(\mu_i)$-a.e..
Let $\cal C$ denote
 the family of all closed subgroups of $T$. We recall that $C$ is
 countable. For $C\in\cal C$, let $I_C=\{\xi\in G/T\mid S_\xi\subset
 C\}$.  Then $I_C$ is a closed subgroup of $G/T$.
%Since $q^{-1}(M)$ centralises $R$, $g_\xi=e$ for all $\xi\in R/T$, we get that
% $R/T$ is contained in $ I_C$ for all $C\in\cal C$.
We have  $\sigma
 (\cup_{C\in\cal C}I_C)=1$, and hence for every index $i$, $q(\mu_i)
 (\cup_{C\in\cal C}I_C)=1$.

Now consider any index $i$. By the above observation
there exists $C\in \cal C$ such that
 $q(\mu_i) (I_C)>0$.  Then $\theta (\mu_i)(\theta(q^{-1}(I_C)))>0$, and
 since $\theta (\mu_i)^{d_i}$ has a density when viewed as a measure on
 $W_i$, this implies that the closed
subgroup $\theta(q-1(I_C)\cap W_i$ has positive Lebesgue
measure in $W_i$, hence it equals $W_i$ and so $\theta
 (\mu_i)(\theta(q^{-1}(I_C))=1$. On the other hand $q^{-1}(I_C)$ is a
 subgroup containing $R=\ker \theta$, and hence  $\theta
 (\mu_i)(\theta(q^{-1}(I_C)))=\mu_i(q^{-1}(I_C))=q(\mu_i)(I_C) $, so
 we have $q(\mu_i)(I_C)=1$.
Let ${\cal C}'$ be the family
 of all $C\in \cal C$ such that $q(\mu_i) (I_C)=1$; by our observations
 above ${\cal C}'$ is nonempty. It can be seen from dimension
 and
 cardinality considerations that every nonempty collection of
 subgroups from $\cal
 C$ has a minimal element. Let $C_0$ be the minimal element in ${\cal
   C}'$.
 Then $q(\mu_i) (I_{C_0})=1$ and $q(\mu_i) (I_C)=0$ for all $C\in \cal C$
 such that $C\subset C_0$ and $C\ne C_0$. It follows that for
 $q(\mu_i)$-almost all $\xi$, $S_\xi=C_0$. Let $E=\{x\in G\mid
S_{xT}=C_0\}$; then we have $\mu_i (E)=1$. Since
for any $\xi\in G/T$ and $\zeta\in R/T$ we have $g_{\xi\zeta}=g_\xi
g_\zeta=g_\xi$ we get that $ER=E$. Also, as
noted above  $\theta (\mu_i)$ is invariant under the action of $\wt
G(\Phi)$ and hence in particular that of
$G(\Psi)$. Therefore for any $y\in G(\Psi)$, we have
$\mu_i (yE y^{-1})=\mu_i (yE y^{-1}R)=\theta (\mu_i) (\theta (yE y^{-1}))=
\theta (\mu_i)(\theta (E))=\mu_i (E)=1$.

Now let $g\in I(\mu)\cap q^{-1}(M)$ and $y\in G(\Psi)$ be given. In
particular $g\in I(\mu_i)$. Also, $q^{-1}(M)$ being
$G(\psi)$-invariant, we get that $ygy^{-1}\in
q^{-1}(M)$ and its action on the fibre $\xi=xT$ is given by
 $$
 (ygy^{-1})x(ygy^{-1})^{-1}x^{-1}=yg(y^{-1}xy)g^{-1}
 (y^{-1}xy)^{-1}y^{-1}=g(y^{-1}xy)g^{-1}(y^{-1}xy)^{-1},$$
the last equality being a consequence of the element being contained
in the center; thus the action is by
translation by $g_{y(\xi)}$, where $y(\xi)=(y^{-1}xy)T$. It follows
that if $\xi \in G/T$ is such that both $\xi$ and $y(\xi)$ are contained in
$E$, then $\mu_\xi$ is invariant under the action of $ygy^{-1}$. Since
$\mu_i (E)=\mu_i(yE y^{-1})=1$ this indeed holds for $q(\mu_i)$-almost all
$\xi$, which shows that $\mu_i$ is invariant under the conjugation
action of $ygy^{-1}$, namely $ygy^{-1}\in I(\mu_i)$. Since this holds
for all $i$ it follows that $ygy^{-1}\in I(\mu)$. This proves the
proposition.
\hfill $\Box$

\bigskip
\noindent {\it Proof of Theorem~\ref{thm:main}}:
As in the proof of  Theorem~\ref{thm:main-tech}, replacing
$\Psi$ by a root cluster $\Psi^*$ contained in it (consisting of
certain fixed powers) we may assume that such
that $\wt G(p(\Psi))/\wt G(p(\mu))$ is connected. When this holds, by
Proposition~\ref{prop:full} the condition as in
  Theorem~\ref{thm:main-tech} is satisfied and hence it follows that
 there exists
a  rational embedding $\{\mu_r\}_{r\in \Q^+}$ of
 $\mu$ such that $\mu_r\in \Psi$ for all $r\in \Q^+$. This proves the
 theorem.
\hfill $\Box$

\bigskip
\noindent{\it Proof of Corollary~\ref{G=NK}}:
Let $T$ be the maximal torus in $N$. Then there exists
a representation $\rho:G\to GL(n,\R)$, for some $n\in \N$, such that
$\ker \rho$ is
contained in the center of $G$ and $(\ker \rho)^0 =T$
(cf. \cite{N}). Now, $[\rho (G),\rho(G)]$ is an almost algebraic
subgroup (cf. \cite{C}) and also $\rho (N)$  is an
algebraic subgroup (cf. \cite{V}), so we get that $[\rho (G),\rho (G)]\rho
(N)$ is almost algebraic. Let $G'$ be the closed subgroup of
$G$ containing $[G,G]N$ such that $G'/\overline{[G,G]N}$ is the
maximal compact subgroup of $G/\overline{[G,G]N}$. Then $G'$
is a closed connected subgroup of $G$ such that $G/G'$ is a vector
group. The latter implies  in particular that $G'$ contains
$\ker \rho$ and hence $\rho (G')$ is a closed subgroup of $\rho (G)$.
From the choice of $G'$ we now get that $\rho (G')/[\rho (G),\rho
(G)]\rho (N)$  is compact, and since $[\rho (G),\rho (G)]\rho
(N)$ is almost algebraic, this implies furthermore that  $\rho (G')$
is algebraic. We can therefore apply  Theorem~\ref{thm:main} to
$G'$, in place of $G$ there.
Now let $\mu \in P(N)$ and suppose that it is infinitely divisible in
$P(G)$. Then $\mu \in P(G')$. Also, since $G/G'$ has
no elements of finite order, all roots of $\mu$ are supported on
$G'$ and so $\mu$ is infinitely divisible on $G'$. Let $q:G'\to G'/T$
be the quotient homomorphism. As $T$ is the maximal torus in
$N$,  $N/T$ has no nontrivial compact
subgroup. Since any idempotent factor is the
Haar measure of a compact subgroup (see \cite{H}) this
shows that $q(\nu)$ has no nontrivial idempotent factor.
Hence by Theorem~\ref{thm:main} $\mu$
is embeddable in $P(G')$, and hence also in $P(G)$.~\hfill $\Box$

\bigskip
For the proof of Corollary~\ref{walnut} we recall some
known results. Let $G$ be a Lie group. A $\mu \in P(G)$ is said to be
{\it strongly root compact} if the set  $\{\rho^k\mid \rho^n=\mu, 1\leq k
\leq n\}$ is relatively compact in $P(G)$, and it is said to be
{\it factor compact} if its set of (two-sided) factors, namely $\lambda
\in P(G)$ for which there exists $\nu\in P(G)$ such that $\mu=\lambda
\nu=\nu\lambda$, is compact. A factor compact measure is strongly
root compact. A strongly root compact measure which is infinitely
divisible  is embeddable (see \cite{H}, \cite{M-cimpa}).

\begin{proposition}\label{strrootcpt}
Let $V=\R^2$ and let $G$ be the semidirect product of $SL(2,\R)$ and
$V$ with respect to the natural linear action. Let $\mu \in P(G)$ be
such that $\mu \notin P(V)$. Then there exists a closed subgroup $G'$
containing $G(\mu)$ such that the following conditions are satisfied:
(i) if $\rho \in P(G)$ and $\rho^{2n}=\mu$ for some $n\in \N$, then
$\rho^2\in P(G')$, and  (ii) $\mu$ is strongly root compact in $G'$.
%Moreover if $L$ is a closed connected subgroup of $G$ containing $V$
% and $\mu \in P(L)$ then $G'$ may be chosen to be contained in $L$.
\end{proposition}

\proof Let $\eta:G\to SL(2,\R)$ be the quotient homomorphism, and let
$\nu=\eta(\mu)$.
Let $A$ be the subgroup of $SL(2,\R)$ consisting of all diagonal matrices
with positive entries, and let $N$ be the subgroup consisting of all upper
triangular unipotent matrices.
From the structure of subgroups of $SL(2,\R)$ we know that at least
one of
the following must hold: (a) $\wt G(\nu)$ contains a conjugate $AN$
in $SL(2,\R)$; (b)~$\wt G(\nu)$ is a nontrivial compact
subgroup;
or (c) $\wt G^0(\nu)$ is a conjugate of $A$ or $N$ in
$SL(2,\R)$; see \cite{M-sl2} for instance, for an idea of the
proof. If~(a) or~(b) holds then
it can be seen that $Z(\mu)$ is a compact subgroup of $G$. As $G$ is
an algebraic
group this implies that $\mu$ is factor compact (cf. \cite{DM-mz}).
Therefore it is strongly root compact, and the assertion in the
proposition holds for $G'=G$.

To prove the assertion when  Condition~(c) holds we may, after
applying a suitable conjugation, assume that  $\wt G^0(\nu)=A$ or $N$.
Then $\wt N(\nu)/\wt G^0(\nu)$ has a normal subgroup $F$ containing
$2$ or $4$ elements such that the corresponding quotient has no
elements of finite order; if  $\wt G^0(\nu)=A$ then $F$ has $4$
elements and the quotient is trivial, while if  $\wt G^0(\nu)=N$ then
$F$ has $2$  elements and the quotient is isomorphic to $\R$. Let $E$
be the subgroup of $\wt N(\nu)$ containing $\wt G^0(\nu)$ and
such that $F=E/\wt G^0(\nu)$. Then $E^0\subset \wt G(\nu)\subset E$.
Also, from the choice of $F$ we see that $\wt N(\nu)/E$ has no
elements of finite order, and therefore it follows that all roots of
$\nu$ are supported
on $E$. Correspondingly, for $\mu$ we have $\eta^{-1}(E^0)\subset \wt G(\mu)
\subset \eta^{-1}(E)$ and all its roots are supported on $\eta^{-1}(E)$.
If $\wt
G(\nu)=N$ we set $G'=\eta^{-1}(N)$. Since $F$ as above has only two
elements in this case, it follows that  $\rho^2\in P(G')$
for all $\rho$ such that $\rho^{2n}=\mu$ for some $n\in \N$. Also,
$N$ is a connected nilpotent Lie group and hence it follows that
$\mu$ is strongly root compact (see \cite{H}, Chapter~3). Thus
the assertion as in the proposition holds in this case.

Next suppose that $\wt G(\nu)\neq N$. Under this condition, if
 $V\cap G(\mu)$ is nonzero, then $Z(\mu)$ is compact and
 hence we are through as before. We may therefore assume that
$V\cap G(\mu)=\{0\}$. Then $V\cap \wt N(\mu)$
is also trivial, since  if $v\in V$ normalizes $\wt G(\mu)$ then
$v$ is a fixed point of the linear action of $\wt G(\nu)$, which
can not hold for any nonzero vector when $\wt
G(\nu)\neq N$. Hence the restriction of
$\eta$ to $\wt N(\mu)$ is an isomorphism of $\wt N(\mu)$ onto $E$. It
therefore
suffices to prove the desired statement for $\mu \in P(E)$ with
all roots contained $E$. If $\wt G(\mu)$ is nonabelian, which is
possible only in the case $\wt G^0(\nu)=A$, then $Z(\mu)$ is compact
and the assertion in the proposition holds as before. If
$\wt G(\mu)$ is
abelian, we choose $G'$ to be the abelian subgroup of index $2$ in
$E$, and see that the requisite conditions are satisfied for this
choice.  This proves the proposition.~\hfill $\Box$

\bigskip
\noindent{\it Proof of Corollary~\ref{walnut}}: Let $G$ be the
semidirect product of $SL(2,\R)$ and $N$, as in the hypothesis, and
$\mu\in P(G)$ be infinitely divisible. Let $T$ be the maximal torus in
$N$.  We shall show that $\mu$ is
embeddable, and that if infinite divisibility holds in $P(H)$, where
$H$ is a closed connected subgroup containing $G(\mu)$, then an embedding
may also be found in $H$; in this respect it suffices to consider
$H$ containing $T$ since for any closed connected subgroup
$H$, $HT$ is a direct product of $H$ and a subtorus of $T$, and
embeddability on the latter implies embeddability on $H$.
If  $\mu\in P(N)$ then Corollary~\ref{G=NK} implies that $\mu$ is
embeddable; we note also that if $\mu\in P(H)$ for a closed connected
subgroup $H$ of $G$, containing $T$, then $H\cap N$ is contained in
the nilradical of $H$ and hence Corollary~\ref{G=NK} may be applied to
$H$ in place of $G$ to conclude that $\mu$ is embeddable in $P(H)$.

Now suppose that $\mu\notin P(N)$.
Let $q:G\to G/T$ be the quotient homomorphisms. Then
$q(\mu) \notin P(N/T)$. Hence by Proposition~\ref{strrootcpt} we
get that there exists a closed subgroup $G'$ of $G/T$ containing
$G(q(\mu))$ such  that for  $\rho \in P(G/T)$ we have $\rho^2\in
P(G')$ whenever $\rho^{2n}=q(\mu)$ for some $n\in \N$, and $q(\mu)$ is
strongly root compact on $G'$. Let $L=q^{-1}(G')$. Then $\mu \in
P(L)$ and  if $n\in \N$ and $\rho \in
P(G)$ is a $2n$th root of $\mu$ then  $\rho^2\in P(L)$. This shows that
$\mu$ is infinitely divisible in $P(L)$. On the other hand since
$q(\mu)$ is strongly root compact on $G'$, and the kernel of $q$
is a compact central subgroup it follows that $\mu$ is strongly
root compact on $L$ (see \cite{M-cimpa}). Together with infinite
divisibility in $P(L)$ this implies that $\mu$ is embeddable in
$P(L)$ (see \cite{H}, Chapter~3, or \cite{M-cimpa}). Hence $\mu$
is embeddable in $P(G)$. Also if $H$ is a closed connected subgroup of
$G$ containing $T$ and such that $\mu$ is infinitely divisible on $H$,
we can choose $L=H\cap q^{-1}(G')$ and arguing as above conclude that
$\mu$ is embeddable in $P(L)\subset P(H)$.  This proves the
Corollary.~\hfill $\Box$

\bigskip
\noindent {\it Proof of Theorem~\ref{noncomm}}. Let $G$ be a connected
Lie group and $\mu \in P(G)$ be an infinitely divisible measure such
that $G(\mu)/\overline{[G(\mu),G(\mu)]}$ is compact. Let $H$ be the
closed subgroup of $G$ containing $[G,G]$ and such that $H/\overline
{[G,G]}$ is the maximal compact subgroup of the connected abelian Lie
group $G/\overline {[G,G]}$. The
condition on $\mu$ implies that $G(\mu)$ is contained in
$H$. Furthermore, since $G/H$ is a vector group we get that all roots
of $\mu$ are supported on $H$, so $\mu$ is infinitely divisible on
$H$. Therefore in proving the theorem, by suitably replacing $G$ we
may assume that $G/\overline {[G,G]}$ is compact. When this condition
is satisfied $G$ admits a representation $p:G\to \wt G$ onto an almost
algebraic group $\wt G$, such that $\ker p$ is contained in the center
of $G$ and $(\ker p)^0$ is compact; this may be seen by arguing as in
the first part of the proof of Corollary~\ref{G=NK}, and noting that
compact connected linear groups are algebraic. We then proceed as in
the proof of Theorem~\ref{thm:main}. In view of the assumption that
$G(\mu)/\overline{[G(\mu),G(\mu)]}$ is compact, the vector space $\wt
V$ as in the earlier sections is trivial, and $\mu$ is in fact
supported on the subgroup $R$ as in~\S~6. For such a measure it
may be seen from the proof of Theorem~\ref{thm:main-tech} that the
condition that $q(\mu)$ has
no idempotent factor is not involved in completing the proof of the
conclusion that every root cluster of $\mu$ contains a rational
embedding of $\mu$. Thus the argument shows that $\mu$ is embeddable
in $P(G)$ under the conditions as above.~\hfill $\Box$

\medskip
\noindent{Acknowledgement}: The first named author would like to thank
the Institut de Recherche Math\'ematique de Rennes (IRMAR),
Universit\'e de Rennes I, France for hospitality while some of this
work was done.

\end{document}